\def\bibart#1#2#3#4#5#6#7
{\bibitem{#1}
#2 (#3), #4, {\it #5} #6, #7}

\def\bibcoll#1#2#3#4#5#6#7#8
{\bibitem{#1}
#2 (#3), #4, in: #5, {\it #6}, #7, p. #8}

\def\bibbook#1#2#3#4#5#6
{\bibitem{#1} #2 (#3), {\it #4}, #5, #6}

\def\bibdiss#1#2#3#4#5#6
{\bibitem{#1}
#2 (#3), #4, #5, #6}

\def\rank{{\rm rank}}

\def\dualrank{{\rm dual rank}}

\def\qed{\hfill\ \rule{2mm}{2mm} }

\documentclass[12pt]{article}

\usepackage{amssymb}

\newtheorem{guess}{Guess}[section]

\newtheorem{define}[guess]{Definition}
\newtheorem{prop}[guess]{Proposition}
\newtheorem{theorem}[guess]{Theorem}

\newtheorem{cor}[guess]{Corollary}

\newtheorem{lem}[guess]{Lemma}

\newtheorem{remark}[guess]{Remark}

\newtheorem{oques}[guess]{Open Question}

\usepackage[pdftitle={Nonextremal Cards}, pdfauthor={Bernd Schroeder}, baseurl={http://www.math.usm.edu/schroeder/}, colorlinks=true, pdfstartview=FitV, linkcolor=blue, citecolor=blue, urlcolor=blue, naturalnames=true, linktocpage=true, bookmarksnumbered=true, pdfpagelayout=SinglePage]{hyperref}

\begin{document}

\bibliographystyle{plain}

\title{
Reconstruction of
the Ranks of the Nonextremal Cards
and of
Ordered Sets with a Minmax Pair of Pseudo-Similar Points
}

\author{
\small Bernd S. W. Schr\"oder \\
\small Department of Mathematics\\
\small The University of Southern Mississippi\\
\small
118 College Avenue, \#5045\\
\small Hattiesburg, MS 39406\\
}

\date{\small \today}

\maketitle

\begin{abstract}

For every
ordered set, we reconstruct the deck
obtained by removal of the
elements of rank $r$ that are neither minimal nor maximal.
Consequently, we also
reconstruct the
deck obtained by removal of the extremal, that is, minimal or maximal,
elements.
Finally, we
reconstruct the
ordered sets with a minmax pair of pseudo-similar points.

{\em 
This preprint has not undergone peer review or any post-submission improvements or corrections. 
The Version of Record of this article is published in Order, and is available online at \url{https://doi.org/10.1007/s11083-024-09684-8}.
}

\end{abstract}

\noindent
{\bf AMS subject classification (2000):} 06A07, 05C60\\
{\bf Key words:}
Reconstruction, ordered set, maximal card, minimal card,
isomorphic cards, rank of the removed element

\section{Introduction}

An {\bf ordered set} consists of a set $P$ and a reflexive,
antisymmetric and transitive relation $\leq $ on $P$,
called the {\bf order relation}.
Subsets of ordered sets inherit the order relation from the
superset.
For a finite ordered set
$P$ and a point $x\in P$, we call the
unlabelled
ordered subset $P\setminus \{ x\} $
a {\bf card} of $P$.
For any finite ordered set $P$, we call the multiset
of unlabelled subsets
${\cal D} _P :=\langle P\setminus \{ x\} :x\in P\rangle $
the {\bf multiset of unlabelled cards} or the {\bf deck} of $P$.
The Reconstruction Problem asks the following question.

\begin{oques}
\label{recprob}

The {\bf Reconstruction Problem.}
Let $P$ and $Q$ be two finite
ordered sets with equal multisets of unlabelled cards
and at least $4$ elements each.
Must $P$ be isomorphic to $Q$?
The {\bf Reconstruction Conjecture} proposes
that the answer is positive.

\end{oques}

Recall that a {\bf minimal} element of an ordered set
$P$ is an element $l$ such that
there is no $p\in P$ with $p<l$.
Dually, a {\bf maximal} element of
$P$ is an element $h$ such that
there is no $p\in P$ with $p>h$.
The set of minimal elements will be denoted $\min (P)$, the
set of maximal elements will be denoted $\max (P)$.
An ordered set $P$ has a {\bf minmax pair $(l,h)$
of pseudo-similar points} iff
$l\in \min (P)$,
$h\in \max (P)$,
and
$P\setminus \{ h\} $ is isomorphic to $P\setminus \{ l\} $.
An element that is maximal or minimal will be called {\bf extremal}.
Moreover, recall that
a {\bf chain} of {\bf length} $k$ is a $(k+1)$-element ordered set
such that any two elements are comparable, and that
the {\bf rank}
of an element
$x$, denoted
${\rm rank}(x)$, is the length
of the
longest chain from a minimal element to $x$.
The largest rank of an element in an ordered set $Q$ is called the
{\bf height} of $Q$, denoted ${\rm height} (Q)$.

Let $\pi $ be a property for points in an ordered set.
We
define
the multiset
${\cal D} _P ^\pi
:=
\langle P\setminus \{ x\} :x\in P, x {\rm \ has \ property \ } \pi \rangle $
to be
the {\bf multiset of unlabelled $\pi $-cards}
or the
{\bf $\pi $-deck} of $P$.
A parameter $\alpha (\cdot )$ of ordered
sets
is called {\bf reconstructible} iff, for any ordered sets $P,Q$,
if ${\cal D} _P = {\cal D} _Q$ then $\alpha (P)=\alpha (Q)$.
Note that $\pi $-decks are parameters.
A class ${\cal K}$ of ordered
sets is called {\bf recognizable} iff, for any ordered sets $P,Q$,
the equality
${\cal D} _P = {\cal D} _Q$
implies that
either both are in ${\cal K}$ or
neither of them is in ${\cal K}$.
An ordered
set $P$ is called {\bf reconstructible} iff, for any ordered set $Q$,
the equality
${\cal D} _P = {\cal D} _Q$
implies that
$P$ is isomorphic to $Q$.

In this paper, we prove the following two theorems.

\begin{theorem}
\label{getnonexcndeck}

Let $P$ be a finite connected ordered set with
at least $4$
elements
and let $r
>0$.
Then the
deck
${\cal N} _P ^r :=\langle P\setminus \{ x\} :
x\in P\setminus \max (P), {\rm rank} (x)=r
\rangle $
of cards obtained by removing
nonmaximal points of
rank $r>0$
is reconstructible.
Consequently, the
deck
$\langle P\setminus \{ x\} :
x\in \max (P)\cup \min (P)
\rangle $
of cards obtained by removing
extremal points
is reconstructible.

\end{theorem}

\begin{theorem}
\label{minmaxpsrecon}

Ordered sets with a minmax pair of pseudo-similar points are
reconstructible.

\end{theorem}

Theorem \ref{minmaxpsrecon} can be seen as the completion of the
answer to
Question 51 of \cite{JXsurvey}, which asks
for a characterization of
the ordered sets (of height $1$) that have a
minmax pair of pseudo-similar points.
A recursive characterization of these sets
was already given in Proposition 3.2 of \cite{bman1}.
We will rely extensively on results from \cite{bman1}, which will be
listed in Section \ref{backgrnd}.
Although, by Corollary 3.5 in \cite{bman1},
reconstructibility would follow easily from the
identification of the maximal deck and the minimal deck,
Theorem \ref{getnonexcndeck} only
identifies
the extremal deck.
Corollary \ref{getsomemincards} indicates that certain
subsets with a minmax pair of pseudo-similar points
and certain chains
are the main obstacles to identifying the maximal deck.
Similarly, we need
insights related to
ordered sets with a minmax pair of pseudo-similar points
in Section \ref{specclassrecon} to
prove Theorem \ref{getnonexcndeck} in Section \ref{nencranks}.
Theorem \ref{getnonexcndeck}
is then used to
prove
Theorem \ref{minmaxpsrecon}
in Sections \ref{conswhenconn} and \ref{isobetwcards}.
Consequences of these results are discussed in Section \ref{conclusionsec}.

\section{Decomposable Ordered Sets}
\label{NTMAcarsec}

In addition to
a needed piece
of the proof of Theorem \ref{getnonexcndeck}
in Proposition \ref{getranksofNTMAcards},
this section
provides further review of
fundamental ideas for later proofs.

For $A,B\subseteq P$, we write $A<B$ iff
every $a\in A$ is strictly smaller than every $b\in B$.
In case one of $A,B$ is a singleton $\{ x\} $, we write
$x<B$ or $A<x$.
Recall that a nonempty
subset $A$ of an ordered set $P$ is
{\bf order-autonomous} iff, for all $p\in P\setminus A$, if there is
an $a\in A$ with $p> a$, then $p>A$;
and,
if there is
an $a\in A$ with $p< a$, then $p<A$.
Singletons and the set $P$ itself trivially are order-autonomous.
An ordered set is called {\bf decomposable} iff it has a nontrivial
order-autonomous subset.

Recall that a {\bf fence from $f_0 $ to $f_n $}
is an ordered set $F=\{  f_0 <f_1 >f_2 <\cdots f_n \} $
or
$F=\{  f_0 >f_1 <f_2 >\cdots f_n \} $
and that an ordered set is called {\bf connected} iff, for any two
elements $a,b\in P$, there is a fence from $a$ to $b$.
Recall that, for two
nonempty ordered sets $A$ and $B$, the
{\bf linear sum}
$A\oplus B$ is the union of $A$ and $B$ with
added comparabilities such that $A<B$,
and that an ordered set is called {\bf coconnected} iff
it cannot be decomposed into a linear sum.
By Theorems 4.3 and 4.7 in \cite{KRtow}, ordered sets that
have a linear sum decomposition or that
are
disconnected
are
reconstructible. Hence, we can assume
that all ordered sets in this paper
are connected and coconnected.

An order-autonomous set is called {\bf maximal} iff it is not
properly contained in a nontrivial order-autonomous set.
Connected and coconnected ordered sets have a canonical decomposition
into a union of pairwise disjoint maximal order-autonomous subsets
$P_t $ with $t\in T$, where $T$ is called the {\bf index set}, ordered by
$t<u$ iff every element of $P_t $ is below every element of
$P_u $.
For every $x\in P$, the index $I(x)\in T$ is the unique element of $T$
such that $x\in P_{I(x)} $.

Although \cite{SchrSetRec}
gives a clear impression that reconstruction of
decomposable ordered sets will follow its own unique path,
to prove Theorem \ref{getnonexcndeck},
we must
include
{\bf NTMA-cards}, that is,
cards
$P\setminus \{ x\} $
for which the removed point $x$ is in
a nontrivial order-autonomous subset of $P$.
Recall that ${\rm Aut} (Q)$ denotes the automorphism group of
an ordered set $Q$.

\begin{lem}
\label{decomplem}

As indicated before Theorem 1.8 in
\cite{SchrSetRec},
all results in
\cite{SchrSetRec} have natural analogues for ordered sets.
In particular, the following hold.
\begin{enumerate}
\item
\label{decomplem1}
(See Theorem 1.8 in \cite{SchrSetRec}. Recognizability
under the assumption $|P|>11$ was proved in \cite{Illedec}.)
Decomposable ordered sets are
(set) recognizable.

\item
\label{decomplem2}
(See Lemma 5.3 in
\cite{SchrSetRec}.)
The NTMA-cards are identifiable in the deck.
In particular, the number of elements in nontrivial order-autonomous
subsets is reconstructible.

\item
\label{decomplem4}
(See Theorem 7.1 in \cite{SchrSetRec}.)
The deck of maximal order-autonomous subsets is
reconstuctible.

\item
\label{decomplem5}
(See Theorem
7.3 in
\cite{SchrSetRec}.)
If at least 3 elements are in nontrivial order-autonomous subsets, then,
for any NTMA-card $P\setminus \{ x\} $,
the ${\rm Aut} \left( T^C \right) $-orbit
of the index set $T^C $ that contains $I(x)$ is identifiable,
the isomorphism type of the
maximal order-autonomous subset $A$ that contains $x$ is
identifiable, as is
(via Theorem 7.1 in \cite{SchrSetRec})
the isomorphism type of $A\setminus \{ x\} $.

\item
\label{decomplem6}
(See Corollary 7.5 in
\cite{SchrSetRec}.)
If the index set $T$ of $P$ is {\bf rigid}
(that is,
the identity is the only automorphism of $T$)
and
at least $3$ elements are in
nontrivial order-autonomous subsets,
then $P$ is reconstructible.
\qed
\end{enumerate}

\end{lem}

Moreover, recall that the number of copies
of a proper subset contained in $P$ is reconstructible.

\begin{prop}
\label{kellylem}

(A Kelly Lemma,
see \cite{KRtow}, Lemma 4.1.)
Let $P$ and $Q$ be two finite ordered sets with $|P|>|Q|$ and $|P|>3$.
Then the number
$s(Q,P)$ of ordered subsets of $P$ that are isomorphic to $Q$
is reconstructible.
\qed

\end{prop}

\begin{define}
\label{neighofA}

Let $P$ be an ordered set and let $A\subseteq P$.
Then we define
$N(A)= \{ p\in P: (\exists a\in A) \ a\leq p {\rm \ or\ } a\geq p\} .$
For $x\in P$, we will write $N(x)$ instead of $N(\{ x\} )$.
When needed, we will
use subscripts
to indicate the surrounding ordered set $P$.

\end{define}

\begin{prop}
\label{getranksofNTMAcards}

Let $P$ be a connected coconnected decomposable ordered set
with at least $3$ elements in nontrivial order-autonomous
subsets, and let
$r\geq 0$.
Then the
deck
${\cal NTMA} _P ^r :=\langle P\setminus \{ x\} :
x\in P, {\rm rank} (x)=r, P\setminus \{ x\} {\rm \ is \ an \ NTMA-card} \}
\rangle $
of unlabelled
NTMA-cards obtained by removing an element of rank $r$
is reconstructible.

\end{prop}

{\bf Proof.}
By part \ref{decomplem2} of Lemma \ref{decomplem}, we can let
$C:=P\setminus \{ x\} $ be an NTMA-card
of $P$.
Let $T^C $ be the index set of its canonical decomposition.
By
part \ref{decomplem5} of Lemma \ref{decomplem},
the isomorphism type of the maximal order-autonomous set $P_{I(x)} $
of $P$
that contains $x$,
the isomorphism type of $P_{I(x)} \setminus \{ x\} $,
as well as the ${\rm Aut} \left( T^C \right) $-orbit $O_x $
that contains $I(x)$
can be identified from the deck.
Consequently, there is a
maximal order-autonomous set
$A\subseteq C$ such that
$I[A]\in O_x $
and $A$ is isomorphic to
$P_{I(x)} \setminus \{ x\} $.

Let $A\subseteq C$ be maximal order-autonomous such that
$I[A]\in O_x $ and $A$ is isomorphic to
$P_{I(x)} \setminus \{ x\} $.
Let $L,U\subseteq C$ be so that
$N_C (A)=L\oplus A\oplus U$.
Then
$N_P \left( P_{I(x)} \right) $
is isomorphic to
$L\oplus P_{I(x)} \oplus U$
iff $C$ contains fewer copies of
$L\oplus P_{I(x)} \oplus U$
than $P$, which is recognizable by Lemma \ref{kellylem}.
Therefore, we can assume that we have reconstructed
sets $L$ and $U$
such that
$N_P \left( P_{I(x)} \right) $
is isomorphic to
$L\oplus P_{I(x)} \oplus U$.

If
the set $L$ of strict lower bounds of
$P_{I(x)} $ is not empty, let
$h$ be the height of
$L$;
otherwise set $h:=-1$.
For every
$s\geq 0$, let $m_s $ be the number of elements
$y\in P_{I(x)} $
whose
rank in
$P_{I(x)} $ is $s$
and
such that
$P_{I(x)} \setminus \{ y\} $
is isomorphic to
$P_{I(x)} \setminus \{ x\} $.
Let $d:=\sum _{s\geq 0} m_s $ and let $t$ be the
total number of cards isomorphic to $P\setminus \{ x\} $.
For every $s\geq 0$,
the
number of
cards isomorphic to
$P\setminus \{ x\} $
in the
deck
${\cal NTMA} _P ^{s+h+1} $
is $m_s $ times ${t\over d} $.
\qed

\section{Sets with Isomorphic Cards}
\label{backgrnd}

This section summarizes the results from
\cite{Schneighdec,bman1}
that we need
for this paper.
Recall that the {\bf components} of an ordered set
are its maximal (with respect to inclusion)
connected subsets.

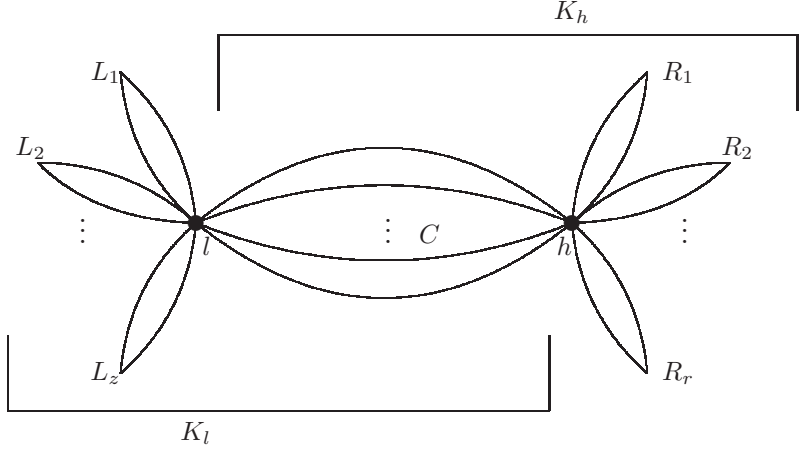
\begin{figure}

\centerline{
\unitlength 1mm 
\linethickness{0.4pt}
\ifx\plotpoint\undefined\newsavebox{\plotpoint}\fi 
\begin{picture}(105,58)(0,0)
\put(25,30){\circle*{2}}
\put(75,30){\circle*{2}}
\qbezier(25,30)(50,40)(75,30)
\qbezier(75,30)(50,20)(25,30)
\qbezier(25,30)(24,42)(15,50)
\qbezier(15,50)(16,38)(25,30)
\qbezier(25,30)(24,18)(15,10)
\qbezier(15,10)(16,22)(25,30)
\qbezier(75,30)(76,42)(85,50)
\qbezier(85,50)(84,38)(75,30)
\qbezier(75,30)(76,18)(85,10)
\qbezier(85,10)(84,22)(75,30)
\put(27,27){\makebox(0,0)[rc]{\footnotesize $l$}}
\put(73,27){\makebox(0,0)[lc]{\footnotesize $h$}}
\put(13,50){\makebox(0,0)[cc]{\footnotesize $L_1 $}}
\put(3,40){\makebox(0,0)[cc]{\footnotesize $L_2 $}}
\put(13,10){\makebox(0,0)[cc]{\footnotesize $L_z $}}
\put(10,30){\makebox(0,0)[cc]{$\vdots $}}
\put(90,30){\makebox(0,0)[cc]{$\vdots $}}
\put(87,50){\makebox(0,0)[lc]{\footnotesize $R_1 $}}
\put(97,40){\makebox(0,0)[cc]{\footnotesize $R_2 $}}
\put(87,10){\makebox(0,0)[lc]{\footnotesize $R_r $}}
\put(50,30){\makebox(0,0)[lc]{$\vdots $ \hspace{.03in} \footnotesize $C$}}
\put(72,15){\line(0,-1){10}}
\put(72,5){\line(-1,0){72}}
\put(0,5){\line(0,1){10}}
\put(28,45){\line(0,1){10}}
\put(28,55){\line(1,0){77}}
\put(105,55){\line(0,-1){10}}
\put(75,58){\makebox(0,0)[cc]{\footnotesize $K_h $}}
\put(25,2){\makebox(0,0)[cc]{\footnotesize $K_l $}}
\qbezier(75,30)(84.1,38)(96,38)
\qbezier(25,30)(15.9,38)(4,38)
\qbezier(96,38)(88.5,30.4)(75,30)
\qbezier(4,38)(11.5,30.4)(25,30)
\qbezier(25,30)(50,50)(75,30)
\qbezier(25,30)(50,10)(75,30)
\end{picture}
}

\caption{Generic picture of a connected ordered set $Q$ with a
minmax pair $(l,h)$ of pseudo-similar points. The view is from
above, with components of $Q\setminus \{ l,h\} $ indicated by ovals.
The sets $K_l $ and $K_h $ are indicated with braces.
}
\label{lhset}

\end{figure}

\begin{define}
\label{lhsetdef}

(See Definition 3.1 of \cite{bman1}.)
Let $Q$ be a finite ordered set
with
a minmax pair $(l,h)$ of
pseudo-similar points.
Throughout
this paper we will use the following notation
(also consider Figure \ref{lhset}).
\begin{enumerate}
\item
\label{lhsetdef1}
$\Phi :Q\setminus \{ h\} \to Q\setminus \{ l\} $
denotes an isomorphism,
\item
$C\subseteq Q$ is the set of all points
$c\in Q\setminus \{ l,h\} $ such that
there is a fence from $c$ to $l$ that does not contain $h$ and
a fence from $c$ to $h$ that does not contain $l$;
the components of $C$ are denoted $C_1 , \ldots , C_c $,
\item
$L\subseteq Q$ is the set of $x\in Q\setminus \{ l,h\} $
such that every fence from $x$ to $h$ contains $l$;
the components of $L$ are denoted $L_1 , \ldots , L_z $,
\item
$R\subseteq Q$ is the set of $x\in Q\setminus \{ l,h\} $
such that every fence from $x$ to $l$ contains $h$;
the components of $R$ are denoted $R_1 , \ldots , R_r $,
\item
$K_l \subseteq Q$ is the component of $Q\setminus \{ h\} $ that contains
$l$,
\item
$K_h \subseteq Q$ is the component of $Q\setminus \{ l\} $ that contains
$h$.
\end{enumerate}
Note that it is possible that $C=\emptyset $ or that $r=0$ or $z=0$.

\end{define}

The following lemma illuminates
the relationship between the various
subsets of a connected ordered set with a minmax pair
$(l,h)$ of pseudo-similar points identified in
Definition \ref{lhsetdef}.

\begin{lem}
\label{componiso}

(See Proposition 3.2 of \cite{bman1}.)
Let $Q$ be a finite connected ordered set with a minmax pair of pseudosimilar
points $(l,h)$.
Then
\begin{enumerate}
\item
$\Phi [K_l ]=K_h $,
\item
$r=z$,
\item
After possibly
renumbering the $L_i $ we have
$\Phi [R_i ]=L_i $ for $i=1, \ldots , z$,
\item
\label{largcomparelh}
$K_l \setminus \{ l\} $ is isomorphic to
$\displaystyle{ K_l \setminus \left\{ \Phi ^{-1} (h) \right\} } $,
that is,
$\displaystyle{ \left( l, \Phi ^{-1} (h)\right) } $
is a minmax pair of pseudo-similar
points in $K_l $.
Similarly, $\displaystyle{ \left( \Phi (l), h\right) } $
is a minmax pair of pseudo-similar
points in $K_h $.
\qed
\end{enumerate}

\end{lem}

The only ways for an ordered set $Q$
to have two elements $l,h$ of
different ranks such that $Q\setminus \{ h\} $ is isomorphic to
$Q\setminus \{ l\} $
are that $(l,h)$
is a minmax pair of pseudo-similar points
or that $Q$ contains an order-autonomous subset that has a minmax
pair of pseudo-similar points.

\begin{prop}
\label{diffrankisdecomp}

(See Proposition 2.3
and
Theorem 2.4 of \cite{bman1}.)
Let $P$ be a finite ordered set.
There are $l,h\in P$ such that
${\rm rank } (l)<{\rm rank } (h)$,
and
$P\setminus \{ l \} $ is isomorphic to $P\setminus \{ h\} $
iff
$P$ contains an
order-autonomous connected
subset $Q$
which contains
a minmax pair $(l,b)$ of pseudo-similar points.
\qed

\end{prop}

Recall that $\downarrow x=\{ p\in P:p\leq x\} $
is the {\bf ideal} of $x$
and that $\uparrow x=\{ p\in P:p\geq x\} $.
The {\bf ideal size sequence} lists all sizes
$|\downarrow x|$ with multiplicity in nondecreasing order.
By Proposition 4.3 in \cite{Schneighdec}, it is reconstructible, and,
by Theorem \ref{keyiso} below,
it uniquely determines
sets with minmax pairs of pseudo-similar points.

\begin{theorem}
\label{keyiso}

(See Theorem 3.4 of \cite{bman1}.)
For $i=1,2$ let $Q_i $ be a finite connected ordered set with a
minmax pair $(l_i , h_i )$ of pseudo-similar points.
If $Q_1 $ and $Q_2 $ have equal ideal size sequences,
then there is an isomorphism $\Psi :Q_1 \to Q_2 $ such that
$\Psi (l_1 ) = l_2 $ and $\Psi (h_1 )=h_2 $.

\end{theorem}

The remaining properties of ordered sets with a minmax pair
of pseudo-similar points proved in \cite{bman1} are best listed
sequentially.

\begin{lem}
\label{propertysequence}

Let $Q$ be a finite connected ordered set with a minmax pair $(l,h)$
of pseudo-similar points
\begin{enumerate}
\item
\label{iterate}
(See Lemmas 2.2 and 3.3 of \cite{bman1}.)
Let
$\Phi :Q\setminus \{ h\} \to Q\setminus \{ l\} $
be an isomorphism.
Then, for all $j\in \{ 0,\ldots , |Q|-1\} $, $\Phi ^j (l)$ is defined,
and,
for
every
$y\in Q$,
there is a unique $j_y \in \{ 0,\ldots , |Q|-1\} $ such that
$\Phi ^{j_y } (l) = y$.
\item
\label{locked}
(See Theorem 4.1 of \cite{bman1}.)
$(l,h)$ is the only minmax pair of pseudo-similar points in $Q$.
\item
\label{stuffaboutlh}
(See Lemma 4.2 of \cite{bman1}.)
No two components of $Q\setminus \{l\} $ are of the same size.
\item
\label{lharerigid}
(See Theorem 4.3 of \cite{bman1}.)
$Q$ and all components of $Q\setminus \{ l\} $ are rigid.
[Note: Parts \ref{stuffaboutlh} and \ref{lharerigid}
imply that $\Phi $ from part \ref{lhsetdef1} of Definition \ref{lhsetdef}
is unique.]
\item
\label{connlhindecomp}
(See Lemma 4.4 of \cite{bman1}.)
If $Q$ is not a chain, then
$\{ l\} $ and $\{ h\} $ are
maximal order-autonomous proper subsets of $Q$.
\item
\label{twocardsonly}
(See Theorem 4.5 of \cite{bman1}.)
If $Q$ is not a chain,
then there is no $a\in Q\setminus \{ l,h\} $ such that
$Q\setminus \{ a\} $ is isomorphic to $Q\setminus \{ l\} $.
\item
\label{dualautlh}
(See Proposition 4.6 of \cite{bman1}.)
$Q$ has a unique dual automorphism $D$.
In particular, $D$ interchanges $l$ and $h$.
\qed
\end{enumerate}

\end{lem}

Finally, we summarize needed results
from
\cite{Schneighdec} in less opaque language.
Recall that $x$ is a {\bf lower cover} of $y$
(and $y$ is an {\bf upper cover} of $x$)
iff $x<y$ and there is no $z\in P$ such that
$x<z<y$.

\begin{lem}
\label{nhoodfacts}

Let
${\cal D} _P =\langle [P\setminus \{ x\} ]:x\in P\rangle $
be the deck of an ordered set $P$. Then the following hold.
\begin{enumerate}
\item
\label{nhoodfacts1}
(See Theorem 3.5 of \cite{Schneighdec}.)
For every $k\geq 0$, the deck of the
unlabelled
neighborhoods of the points of rank $k$
in $P$ is reconstructible.
\item
\label{nhoodfacts2}
(See Proposition 4.3 in \cite{Schneighdec}.)
The deck of the
unlabelled
ideals of $P$ is reconstructible.
\item
\label{nhoodfacts3}
(See Proposition 4.15 of \cite{Schneighdec}.)
The deck of maximal cards
whose
removed
maximal element has
a unique lower cover is reconstructible.
\qed
\end{enumerate}

\end{lem}

\section{A Class of Reconstructible Ordered Sets}
\label{specclassrecon}

As was mentioned in the introduction, the first step for our proofs is the
reconstruction of some special classes of ordered sets that
have proper subsets with
a minmax pair of pseudo-similar points, see
Lemmas \ref{getP-lLemma} and \ref{Qlinear} below.

\begin{lem}
\label{smallcardscor}

(See Corollary 6.4 in \cite{KRtow}.)
Let $P$ be a connected finite ordered set with
at least $4$ elements.
Then it is possible to identify an unlabelled
minimal card $Z=P\setminus \{ z\} $
on which the
minimal elements of $P\setminus \{ z\} $ that are not minimal in $P$
are identified.
Consequently $P\setminus \min (P)$ can be identified on
$P\setminus \{ z\} $ and hence $P\setminus \min (P)$ is reconstructible.
\qed

\end{lem}

\begin{lem}
\label{identcardinupperlevels}

Let $P$ be a connected finite ordered set
with at least $4$ elements
such that there is a $y\in P\setminus \min (P)$
such that the unlabelled card
$Y:=P\setminus \{ y\} $ can be identified,
such that
$Y\setminus \min (P)$ is rigid,
and such that
$y$ can be identified on a card
$Z:=P\setminus \{ z\} $ as in Lemma \ref{smallcardscor}.
Then $P$ is reconstructible.

\end{lem}

{\bf Proof.}
Let $Q:=P\setminus \min (P)$.
On $Z$,
by Lemma \ref{smallcardscor},
we can identify $Q_Z :=Z\setminus \min (P)=P\setminus \min (P)$
and, by assumption,
we can identify $y\in Q_Z $.
On $Y$, because $y$ is not minimal, we can identify
$\min (P)=\min (Y)$, and hence we can identify
$Q\setminus \{ y\} =Y\setminus \min (P)$.
Now let
$\Psi :Y\setminus \min (P)\to Q_Z \setminus \{ y\} $
be
an
isomorphism between these sets.
Because
$Q \setminus \{ y\} =Y\setminus \min (P)= Q_Z \setminus \{ y\} $
is rigid,
$\Psi $ is unique and,
although both cards are unlabelled,
it must be the identity.

Decompose
$\min (P)\subset Y$
into
pairwise disjoint
maximal $Y$-order-autonomous antichains
$A_1 , \ldots , A_m $.
Decompose
$\min (P)\setminus \{ z\} \subset Z$
into pairwise disjoint
subsets
$B_1 , \ldots , B_m $
and
$C_1 , \ldots , C_m $
that are empty or
$Z$-order-autonomous antichains
such that, for all $j\in \{ 1, \ldots , m\} $, we have that
$y\not> B_j $ unless $B_j =\emptyset $,
$y>C_j $, and, if
$
B_j
\not= \emptyset $ and
$
C_j \not= \emptyset $, then
$
\{ u\in Z: u>C_j \}
\setminus \{ y\}
=
\{ u\in Z: u>B_j \}
$.
Because
$\Psi $ must be the identity, the enumerations can be
chosen
such that,
for all $j\in \{ 1, \ldots , m\} $,
if $B_j \not= \emptyset $, then
$
\Psi [\{ x\in Y:x>A_j \} ]
=
\{ u\in Z: u>B_j \}
$
and, if $C_j \not= \emptyset $, then
$
\Psi [\{ x\in Y:x>A_j \} ]
=
\{ u\in Z: u>C_j \} \setminus \{ y\}
$.
Moreover, again
because
$\Psi $ must be the identity, the enumerations can be
chosen
such that,
for all $j\in \{ 1, \ldots , m-1\} $,
we have that
$|A_j |=|B_j |+|C_j |$,
and such that
$|A_m |=|B_m |+|C_m |-1$.

For such an enumeration, we conclude that
$z\in A_m $.
Now, if
$|\uparrow z\setminus \{ z\} |=|\Psi [\{ x\in Y:x>A_m \} ]|$, then
$\uparrow z\setminus \{ z\}
=\Psi [\{ x\in Y:x>A_m \} ]$.
Otherwise, we have
$|\uparrow z\setminus \{ z\} |=|\Psi [\{ x\in Y:x>A_m \} ] |+1$,
and then
$\uparrow z\setminus \{ z\} =\Psi [\{ x\in Y:x>A_m \} ] \cup \{ y\} $.
Because
$\uparrow z\setminus \{ z\} $ is thus
identified on $Z$, we have reconstructed $P$.
\qed

\begin{lem}
\label{getP-lLemma}

Let $P$ be a connected ordered set
with at least $4$ elements
such that $Q:=P\setminus \min (P)$ is
a connected ordered set with a minmax pair $(l,h)$
of pseudo-similar points.
Then $P$ is reconstructible.

\end{lem}

{\bf Proof.}
If $Q$ is a chain, then $P$ has a largest element and is hence
reconstructible. Therefore, we can assume that
$Q$ is not a chain.
Hence, by part \ref{locked} of Lemma \ref{propertysequence},
$l$ is the unique minimal element of
$P\setminus \min (P)$ whose removal
produces a card that is isomorphic to a maximal card of $P\setminus \min (P)$.
That is,
$l$ can be identified on a card
$Z:=P\setminus \{ z\} $ as in Lemma \ref{smallcardscor}.
Let $Y:=P\setminus \{ l\} $.
By parts \ref{stuffaboutlh} and \ref{lharerigid}
of Lemma \ref{propertysequence},
$Y\setminus \min (P)=Q\setminus \{ l\} $
is rigid.
Throughout this proof, it will
be our goal to
prove that $P\setminus \{ l\} $
can be identified in the deck.
Via Lemma \ref{identcardinupperlevels},
this implies that $P$ is reconstructible.

{\em Claim 1.
We can identify
$P\setminus \{ l\} $
or
we can identify
two cards $C$ and $K$ such that
$\{ C,K\} =\{ P\setminus \{ l\} , P\setminus \{ h\} \} $.}
Clearly, the
cards
$X$ such that $X\setminus \min (X)$
is isomorphic to
$Q\setminus \{ l\} $,
are identifiable
in the deck.
If there are exactly two such cards,
then
we have identified
$\{ P\setminus \{ l\} , P\setminus \{ h\} \} $.
Otherwise,
by part \ref{twocardsonly} of Lemma \ref{propertysequence},
the point $l$ must have a unique lower cover
$m$, and then
$P\setminus \{ l\} $,
$P\setminus \{ h\} $,
$P\setminus \{ m\} $,
are the only cards $X$ of $P$
such that $X\setminus \min (X)$ is
isomorphic to
$Q\setminus \{ l\} $.

For $x\in \{ l,h,m\} $, we let
$u_x $ be the sum of the
number of upper bounds
of the
$P\setminus \{ x\} $-minimal elements on
$P\setminus \{ x\} $.
Because $l$ is above exactly one minimal
element of $P$,
$h$ is above at least one minimal element of $P$, and
$|\uparrow m|\geq |\uparrow l|+1$,
we have $u_l\leq u_h , u_m $.
In the identifiable case $u_l< u_h, u_m $, we can identify
$P\setminus \{ l\} $.
If $u_l = u_m $, then
$l$ is the unique upper cover of $m$,
which is identifiable by the dual of
part \ref{nhoodfacts3} of Lemma \ref{nhoodfacts}.
In this case, $\{ l,m\} $
is order-autonomous,
and hence
$P\setminus \{ l\} $ is isomorphic to
$P\setminus \{ m\} $.
Independent of whether $P\setminus \{ h\} $ is isomorphic to the other two
cards or not, in this case, we can therefore
identify two unmarked cards
$C$ and $K$ such that
$\{ C,K\} =\{ P\setminus \{ l\} , P\setminus \{ h\} \} $:
Pick one of the two (or three) isomorphic cards as $C$ and
pick the nonisomorphic card (or another isomorphic one) as $K$.
This only leaves the, now also identifiable, case
$u_l =u_h <u_m $, in which case we can
identify
$\{ C,K\} =\{ P\setminus \{ l\} , P\setminus \{ h\} \} $
as the set of cards $P\setminus \{ x\} $ with the smallest (and equal) values
of $u_x $.
This proves
{\em Claim 1}.

{\em Claim 2.
There are two distinct elements $x,y\in
P\setminus \min (P)$ such that
$
|\downarrow x\cap \min (P)|
\not=
|\downarrow y\cap \min (P)|
$ or $P$ is reconstructible.}
Consider the case that,
for all
$x,y\in
P\setminus \min (P)$, we have
$
|\downarrow x\cap \min (P)|
=
|\downarrow y\cap \min (P)|
$.
Then, for any two
$x,y\in
P\setminus \min (P)$ such that $x<y$, we have
$
\downarrow x\cap \min (P)
\subseteq
\downarrow y\cap \min (P)
$, and,
because
$
|\downarrow x\cap \min (P)|
=
|\downarrow y\cap \min (P)|
$,
we obtain
$
\downarrow x\cap \min (P)
=
\downarrow y\cap \min (P)
$.
Because $P\setminus \min (P)$ is connected, we obtain that, for any two
$x,y\in P\setminus \min (P)$, we have
$
\downarrow x\cap \min (P)
=
\downarrow y\cap \min (P)
$.
This implies that $\min (P)<P\setminus \min (P)$,
which implies that $P=
\min (P)\oplus (P\setminus \min (P))$ is reconstructible.
We have proved
{\em Claim 2}.

With $\{ C,K\} $ from {\em Claim 1},
let $N:=|C\setminus \min (C)|=|K\setminus \min (K)|=|Q|-1$.
Let
$\Psi _l :C\setminus \min (C)\to Q\setminus \{ h\} $
and
$\Psi _h :C\setminus \min (C)\to Q\setminus \{ l\} $
be the, by parts
\ref{stuffaboutlh},
\ref{lharerigid}
and \ref{twocardsonly} of Lemma \ref{propertysequence} only two possible,
isomorphisms from $C\setminus \min (C)$ to cards of $Q$.
Now $\Psi _h \Psi _l ^{-1} =\Phi $ and,
by part \ref{iterate} of Lemma \ref{propertysequence}, for all
$j\in \{ 0,\ldots , N\} $, we have that
$
\Psi _l ^{-1} \left( \Phi ^j (l)\right)
=
\Psi _h ^{-1} \left( \Phi ^{j+1} (l)\right) $.
Therefore,
we can
enumerate the elements of
$C\setminus \min (C)$ as
$x_1 , \ldots , x_N $
such that $i<j$ implies that
$x_j =\Phi ^{j-i} (x_i )$.
Similarly, we can
enumerate the elements of
$K\setminus \min (K)$ as
$y_1 , \ldots , y_N $
such that $i<j$ implies that
$y_j =\Phi ^{j-i} (y_i )$.
For every $j\in \{ 0, \ldots , N\} $, we set
$
c_j :=|\downarrow x_j \cap \min (P)|
$
and
$
k_j :=|\downarrow y_j \cap \min (P)|
$.

For $j=0, \ldots , N$, let
$q_j :=|\downarrow \Phi ^j (l)\cap \min (P)|$.
Although the ordered sequence
$q_0 , \ldots , q_N $
is not reconstructible,
by part \ref{nhoodfacts2} of Lemma \ref{nhoodfacts},
the
unordered multiset of the
$q_j $ is reconstructible.
Hence, for
$C=P\setminus \{ a\} $, the missing
size
$e_1 :=|\downarrow a\cap \min (P)|$
can be identified, and,
for
$K=P\setminus \{ b\} $, the missing
size
$e_2:= |\downarrow b\cap \min (P)|$
can be identified.
Moreover, because $\{ a,b\} =\{ l,h\} $, we know that
$\{ e_1 , e_2 \}
=
\{ |\downarrow l \cap \min (P)| , |\downarrow h \cap \min (P)| \}
$.

One of the sequences
$(e_1 , c_1 , \ldots , c_N )$
and
$(c_1 , \ldots , c_N , e_1 )$
is equal to
$(q_0 , q_1 , \ldots , q_N )$,
and
one of the sequences
$(e_2 , k_1 , \ldots , k_N )$
and
$(k_1 , \ldots , k_N , e_2 )$
is equal to
$(q_0 , q_1 , \ldots , q_N )$.
Because $\{ e_1 , e_2 \}
=
\{ |\downarrow l \cap \min (P)| , |\downarrow h \cap \min (P)| \}
$, among the four sequences
$(e_1 , c_1 , \ldots , c_N )$;
$(c_1 , \ldots , c_N , e_1 )$;
$(e_2 , k_1 , \ldots , k_N )$
and
$(k_1 , \ldots , k_N , e_2 )$,
with $\{ i,j\} =\{ 1,2\} $,
there is a sequence that starts with
$e_i $ that is equal to a sequence that ends with $e_j $.
Without loss of generality, we can assume that
$(c_1 , \ldots , c_N , e_1 )
=
(e_2 , k_1 , \ldots , k_N )$.
Now, if
$(e_1 , c_1 , \ldots , c_N )
\not=
(k_1 , \ldots , k_N , e_2 )
$,
then
$C=P\setminus \{ h\} $ and $K=P\setminus \{ l\} $
has been identified.

This leaves the case that
$(c_1 , \ldots , c_N , e_1 )
=
(e_2 , k_1 , \ldots , k_N )$
and
$(e_1 , c_1 , \ldots , c_N )
=
(k_1 , \ldots , k_N , e_2 )
$.
Then
$e_1 =k_1 =c_2 =k_3 =c_4 =k_5 =\cdots $
and
$e_2 =c_1 =k_2=c_3=k_4 =\cdots $.
Independent of whether
$(e_1, c_1 , \ldots , c_N )
=(q_0 , \ldots , q_N )$
or
$(e_2, k_1 , \ldots , k_N )
=(q_0 , \ldots , q_N )$,
we obtain that, in this case, all
the even indexed
$q_{2j } $ are equal and all
the odd indexed
$q_{2i +1} $ are equal.

In case there are a $j>0$ and an $i\geq 0$ such that
$\Phi ^{2j} (l)<\Phi ^{2i+1} (l)$, we have
$q_{2j} \leq q_{2i+1} $ and, via
$\Phi ^{2j-1} (l)<\Phi ^{2i} (l)$, we have
$q_{2j-1} \leq q_{2i} $.
Consequently, all $q_j $ are equal and, by {\em Claim 2},
$P$ is reconstructible.
Similarly, if there are a $j\geq 0$ and an $i\geq 0$
such that $\Phi ^{2i+1} (l)\not= h $
and $\Phi ^{2j} (l)<\Phi ^{2i+1} (l)$,
then $P$ is reconstructible.

This leaves the case that
$\Phi ^{2j} (l)<\Phi ^{2i+1} (l)$
implies $\Phi ^{2j} (l)=l$
and
$\Phi ^{2i+1} (l)=h$.
In case
$
|\downarrow l\cap \min (P)|
=
|\downarrow h\cap \min (P)|
$, by {\em Claim 2},
$P$ is reconstructible.

We are left with the case that
$\Phi ^{2j} (l)<\Phi ^{2i+1} (l)$
implies $\Phi ^{2j} (l)=l$
and
$\Phi ^{2i+1} (l)=h$
and
$
|\downarrow l\cap \min (P)|
<
|\downarrow h\cap \min (P)|
$.
In particular
$
\downarrow l\cap \min (P)
$
is strictly contained in
$
\downarrow h\cap \min (P)
$.
Because $Q$
is connected, any two even
(odd) indexed $q_j $ are connected via
a fence consisting of even (odd) indexed $q_i $.
Consequently, if $j,k$ are both even (odd), then
$\downarrow q_j \cap \min (P) = \downarrow q_k \cap \min (P)$,
and then every minimal element below
$l$ is below $P\setminus \min (P)$.

Let $m\in \min (P)$
such that $m<h$ and $m\not< l$.
Then $P\setminus \{ m\} $ has
fewer minimal elements than $P$.
Hence
$P\setminus \{ m\} $ can be identified as a minimal card
on which we can identify $P\setminus \min (P)$
(and hence $l$ and $h$),
on which $l$ has as many minimal lower bounds as in $P$,
and
on which $h$ has fewer minimal lower bounds than in $P$.
Because
$l$ has as many minimal lower bounds as in $P$
and
$h$ has fewer minimal lower
bounds in $P\setminus \{ m\} $ than
in $P$, we identify that we must have
$m<h$ and $m\not<l$.
Thus, in this last case,
we reconstruct $P$ by adding $m$ to $P\setminus \{ m\} $
as a strict lower bound for all
$\Phi ^{2i+1} (l)$.
\qed

\begin{lem}
\label{Qlinear}

Let $P$ be a connected ordered set
with at least $4$ elements
such that $Q:=P\setminus \min (P)$ is
a linear
sum $B\oplus A\oplus B'$, where
$A$ is a connected ordered set with a minmax pair $(l,h)$ of pseudo-similar points
such that $A$ is not a chain,
$B$ and $B'$ could each be empty
and each minimal element of $P$ is below some element of $A$.
Then $P$ is reconstructible.

\end{lem}

{\bf Proof.}
We first prove recognizability.
Let $P$ be an ordered set.
By Lemma \ref{smallcardscor}, the set $Q=P\setminus \min (P)$ is
reconstructible.
If $Q$ does not have a linear sum decomposition as indicated,
then $P$ is not as indicated in the lemma.
Otherwise,
let $A$ be the unique linear summand of $Q$ that
is connected, that is not a chain,
that contains a minmax pair $(l,h)$ of pseudo-similar points,
and whose set of upper bounds $B'$ does {\em not} contain
another linear summand that is connected,
contains a minmax pair of pseudo-similar points,
and that is not a chain.

Let $m\in P$ be minimal.
If there is an $a\in A$ such that $m<a$, then $a<B'$ implies
$|\uparrow m|\geq |\uparrow a|+1\geq |B'|+2>|B'|+1$.
If there is no $a\in A$ such that $m<a$, then
$\uparrow m\subseteq \{ m\} \cup B'$ implies
$|\uparrow m|\leq |B'|+1$.
Thus $P$ is an ordered set as indicated iff $Q=P\setminus \min (P)$ has the requisite structure
and, for all minimal elements of $P$ we have $|\uparrow m|>|B'|+1$.
Because the $|\uparrow m|$ are reconstructible
via part \ref{nhoodfacts1} of Lemma \ref{nhoodfacts},
ordered sets as indicated
are recognizable.

In case $B'\not= \emptyset $ (which is recognizable from $Q$),
the set
$P$ is the linear
sum
$(P\setminus B')\oplus B'$ with both summands not empty, which means $P$
is reconstructible.
For the remainder of this proof, we can assume that
$B'=\emptyset $.
In case $B=B'=\emptyset $
(which is recognizable from $Q$), $P$ is reconstructible by
Lemma \ref{getP-lLemma}.
Thus, for the remainder of this proof, we can assume that
$B'=\emptyset $ and $B\not= \emptyset $.

If all minimal elements are lower bounds of
$A$, then $P$ is the linear
sum
$(P\setminus A)\oplus A$ and hence
$P$ is reconstructible.
Thus we can assume that
$P$ has a minimal element that is not a
lower bound of $A$.

All
elements of $B$
are lower bounds of $A$.
Therefore the minimal elements that are not lower bounds of $A$
are not
comparable to any element of $B$.
Because $P$ is connected, for all these minimal elements
$m$, the set $\uparrow m\setminus \{ m\} \not= \emptyset $
is a strict subset of
$A$.
Because no element of $B$ is minimal, there must be a minimal
element that is a lower bound of an element of $B$, which means that the set
$L=\{ p\in P: p\leq A\} $ has at least two elements.
Because
$B\subseteq L$ and the minimal elements that are
not lower bounds of $A$ are not comparable to any element of $B$,
we have that $L$ is order-autonomous in $P$.
Because $|L|\geq 2$, $P$ is decomposable, which,
by
part \ref{decomplem1} of Lemma \ref{decomplem},
is
recognizable.

By
part \ref{decomplem4} of Lemma \ref{decomplem},
the isomorphism types of the maximal order-autonomous subsets of $P$
are reconstructible.
The NTMA-cards $C=P\setminus \{ v\} $
obtained by removing an element $v\in L$ are the
NTMA-cards
such that $C\setminus \min (C)$ is a linear
sum
$\tilde{B}\oplus \tilde{A}$
with $\tilde{A}$ isomorphic to $A$
such that, for all
$c\in C$, there is an $a\in \tilde{A} $ with $a\geq c$
and such that $\tilde{A}$ has as few lower bounds as possible.
In particular, $A$ is identifiable
(as the set $\tilde{A}$) on these NTMA-cards.
Fix one such NTMA-card $P\setminus \{ v\} $.
The set $L\setminus \{ v\} $ is identifiable as the set of lower bounds
of $A$.
Replacing the lower bounds $L\setminus \{ v\} $
of $A$ on $P\setminus \{ v\} $ with an ordered set isomorphic to $L$
as the new set of lower bounds of $A$ (and adding no further comparabilities)
provides an ordered set isomorphic to $P$.
Hence
$P$ is
reconstructible.
\qed

\section{
Proof of Theorem \ref{getnonexcndeck}
}
\label{nencranks}

\begin{lem}
\label{getranknonex1}

Let $P$ be a connected ordered set
with at least $4$ elements
and let $P\setminus \{ x\} $ be a non-NTMA
card of $P$ with the same number of minimal
elements as $P$.
Then
we can reconstruct an
$r\in \{ 1, \ldots , {\rm rank} _P (x)\} $ such that
$\rank (x)\in \{ 0,r\} $.

\end{lem}

{\bf Proof.}
Let
$U_x :=\left( P\setminus \{ x\} \right) \setminus
\min \left( P\setminus \{ x\} \right) $
and let
$Q:=P\setminus \min (P)$.
Clearly, if $x$ is not minimal, then
$U_x = Q\setminus \{ x\} $, whereas, if
$x$ is minimal, then
$U_x = Q\setminus \{ m\} $ for a $Q$-minimal element
$m\not= x$.
Equally clearly,
if all elements $q\in Q$ such that
$Q\setminus \{ q\} $ is isomorphic to $U_x $ have the
same $Q$-rank
$r_Q $,
then
$r:=r_Q +1$ leads to
${\rm rank} _P (x)\in \{ 0,
r\} $.

We are left with the case that there are
elements $q_1 , q_2 \in Q$
that have different
$Q$-ranks and
such that each
$Q\setminus \{ q_i \} $ is isomorphic to $U_x $.
Clearly, this is identifiable, and,
by
Proposition \ref{diffrankisdecomp},
there are two possibilities:
Either
$U_x $ is isomorphic to a card of $Q$ that
was obtained by removal of the bottom element of a
$Q$-order-autonomous chain $C$ with $2$ or more elements, or,
$U_x $ is isomorphic to a card of $Q$ that
was obtained by removal of the
element $l$ of a
connected $Q$-order-autonomous subset $A$ with
a minmax pair
$(l,h)$
of pseudo-similar points that is not a chain.
In either of these identifiable cases,
it is possible that $Q=C$ or $Q=A$, respectively.

{\em Case 1:
$U_x $ is isomorphic to a card of $Q$ that
was obtained by removal of the bottom element of a
$Q$-order-autonomous chain
$C=\{ c_1 <c_2 <\cdots <c_n \} $
with $n\geq 2$ elements.
}
In this case,
the $Q$-rank
$t$ of $c_n $ is identifiable.
Let
$e<t$ be the largest number such that,
for all $j\leq e$, we have
$\displaystyle{
\sum _{{\rm rank} _{P\setminus \{ x\} } (z)=j} |\downarrow _P z|
=
\sum _{{\rm rank} _{P} (z)=j} |\downarrow _P z| } $
and set $r:=t+1$.
Clearly, if $x$ is minimal, then ${\rm rank} (x)\in \{ 0, r\} $.
In case $x$ is not minimal, there is a
$k\in \{ 1, \ldots , n\} $
such that $x=c_k $.
Therefore, for all $j<{\rm rank} _P (x)$, we have
$\displaystyle{
\sum _{{\rm rank} _{P\setminus \{ x\} } (z)=j} |\downarrow _P z|
=
\sum _{{\rm rank} _{P} (z)=j} |\downarrow _P z| } $.
Hence, in particular, if $k=n$, we have that ${\rm rank} _P (x)=r$.
In case $k<n$, because $P\setminus \{ x\} $ is not an NTMA-card,
$\{ c_k , c_{k+1} \} $ is not $P$-order-autonomous.
Because
$\{ c_k , c_{k+1} \} $ is $Q$-order-autonomous, there must be a
$P$-minimal element $z$ such that $z<c_{k+1} $ and $z\not< c_k $.
Because the $P\setminus \{ x\} =P\setminus \{ c_k \} $-rank
$r_{k+1} ^x $
of
$c_{k+1} $ is one less than its $P$-rank
and hence equal to the $P$-rank of $c_k $, we obtain that
$\displaystyle{
\sum _{{\rm rank} _{P\setminus \{ x\} } (z)=r_{k+1} ^x } |\downarrow _P z|
>
\sum _{{\rm rank} _{P} (z)=r_{k+1} ^x } |\downarrow _P z| } $.
Hence, in this case, too,
${\rm rank} _P (x)={\rm rank} _P (c_k )=r_{k+1} ^x =r$.
This concludes {\em Case 1}.

{\em Case 2:
$U_x $ is isomorphic to a card of $Q$ that
was obtained by removal of the element $l$ of a
connected $Q$-order-autonomous subset $A$ with
a minmax pair
$(l,h)$
of pseudo-similar points that is not a chain.
}
Let $N(A)$
be
the
set
from Definition \ref{neighofA}.
We first
consider the
case $N(A)=P$.
In this case, all elements of $Q$ are comparable to some element of
$A$.
Because $A$ is $Q$-order-autonomous,
$Q$ is a linear
sum
$B\oplus A\oplus B'$, where $B,B'$ or both could be empty, and
every minimal element of $P$ is below some element of $A$.
Hence, in this case, by Lemma \ref{Qlinear},
$P$ is reconstructible.
We are left to consider the case that $N(A)\not= P$.

From $Q$ and $U_x $, we can
identify the isomorphism type of $A$ and
the isomorphism type of $N(A)\setminus \min (N(A))$
(which is equal to $N_Q (A)$).
Because $N(A)\not= P$,
by Proposition \ref{kellylem}, we can identify
a largest possible ordered set
$V$ such that $P$ contains more copies of $V$ than $P\setminus \{ x\} $,
such that there is an isomorphism
$\Psi : V\setminus \min (V)\to N(A)\setminus \min (N(A))$
and such that every element of $V$ is comparable to
an element of $\Psi ^{-1} [A]$.
If $x$ is not minimal, then $x\in \{ l,h\} \subseteq A$,
$V$ must be isomorphic to $N(A)$
and the isomorphism is
a natural extension of $\Psi $ to include $\min (V)$.

{\em Claim.
If
$V\setminus \left\{ \Psi ^{-1} (l)\right\} $ is isomorphic to
$V\setminus \left\{ \Psi ^{-1} (h)\right\} $, then
$x$ is minimal.}
Suppose, for a contradiction, that
$V\setminus \left\{ \Psi ^{-1} (l)\right\} $
is isomorphic to $V\setminus \left\{ \Psi ^{-1} (h)\right\} $
and $x$ is not minimal.
Then $x\in \{ l,h\} $ and
$N(A)\setminus \{ l\} $ is isomorphic to
$N(A)\setminus \{ h\} $.
By Proposition \ref{diffrankisdecomp},
$l$ is contained in a nontrivial $N(A)$-order-autonomous
connected ordered set $M$ which,
because
${\rm rank}_P (h)>{\rm rank}_P (l)$, contains a
minmax pair $(l,k)$ of pseudo-similar points with $k\not=l$.
Because $l$ is minimal in $M$,
the smallest $N(A)$-rank of any element of $M$ is ${\rm rank}_{N(A)} (l)$.

We claim that $M\cap A$ is $N(A)$-order-autonomous.
Because the intersection of order-autonomous sets is
order-autonomous, $M\cap A$ is order-autonomous in
$N(A)\setminus \min (N(A))$.
To finish the proof of $N(A)$-order-autonomy, let
$z\in N(A)\setminus (M\cap A)$ be a strict
$N(A)$-minimal
lower bound of an element of
$M\cap A$.
Then
${\rm rank}_{N(A)} (z)=0\leq
{\rm rank}_{N(A)} (l)-1$, which means in particular that
$z\not\in M$.
Now, because $M$ is $N(A)$-order-autonomous,
$z$ is a lower bound of $M$ and hence it is a lower bound of
$M\cap A$.

Because $M\cap A$ is $N(A)$-order-autonomous, it is $A$-order-autonomous and
it contains $l$.
Moreover, because $M$ is connected, not a singleton and contains no elements
of rank less than ${\rm rank} _N (l)$,
$M$ contains at least one upper cover of $l$.
However, all upper covers of $l$ are in $A$, so
$|M\cap A|>1$.
By part \ref{connlhindecomp} of Lemma \ref{propertysequence}
this means
that $M\cap A=A$.
Hence $A$ is $N(A)$-order-autonomous.
However, $N(A)$ contains all points that are comparable to some
point in $A$, so $A$ is $P$-order-autonomous.
Consequently,
$P\setminus \{ x\} $ is an NTMA-card,
a contradiction to our hypothesis.
This establishes the
{\em Claim}.

Because the rank of a minimal element is contained in any
set
$\{ 0,r\} $,
we are left with the case that
$V\setminus \left\{ \Psi ^{-1} (l)\right\} $
is not isomorphic to
$V\setminus \left\{ \Psi ^{-1} (h)\right\} $.

In case
$P\setminus \{ x\} $ contains fewer copies of
$V\setminus \left\{ \Psi ^{-1} (h)\right\} $
than $P$ and the same number of copies
of $V\setminus \left\{ \Psi ^{-1} (l)\right\} $,
we set $r:={\rm rank} _Q (l)+1$, otherwise, we
set $r:={\rm rank} _Q (h)+1$.

If $x$ is not minimal,
then
$x\in \{ l,h\} $.
Now $x=l$ iff $P\setminus \{ x\} $ contains fewer copies of
$N(A)\setminus \{ h\} $
(which is isomorphic to $V\setminus \left\{ \Psi ^{-1} (h)\right\} $)
than $P$ and
the same number of copies
of $N(A)\setminus \{ l\} $
(which is isomorphic to $V\setminus \left\{ \Psi ^{-1} (l)\right\} $)
and in this case ${\rm rank} _P (x)=r$.
Otherwise
$x=h$ and again ${\rm rank} _P (x)=r$.
\qed

\begin{theorem}
\label{nonextcardthmnew}

Let $P$ be a connected ordered set with $\geq 4$ elements.
The nonextremal non-NTMA cards are identifiable
in the deck. Moreover, for each such card, the
rank of the removed element is reconstructible.

\end{theorem}

{\bf Proof.}
Let $P\setminus \{ x\} $ be a non-NTMA card of
$P$.
If the number of
$P\setminus \{ x\} $-maximal elements
does not equal the number of
$P$-maximal elements,
then
$P\setminus \{ x\} $ is identified as
maximal.
Similarly,
if the number of
$P\setminus \{ x\} $-minimal elements
does not equal the number of
$P$-minimal elements,
then
$P\setminus \{ x\} $ is identified as
minimal.
For the remainder, we can assume that
$P\setminus \{ x\} $
has as many minimal elements as $P$ and as many
maximal elements as $P$.

By Proposition \ref{kellylem},
we can reconstruct the length $c>0$ of the
longest chain in $P$ that contains $x$
as
the unique
largest $c>0$ such that $P$ contains
more chains of length $c$ than $P\setminus \{ x\} $.
Recall that the
{\bf dual rank}
of an element
$x$,
denoted
${\rm dualrank}(x)$,
is the length of the
longest chain from $x$ to a maximal element.
By Lemma \ref{getranknonex1} and its dual, we can
reconstruct
$r,d\in \{ 1, \ldots , c\} $
such that
$\rank (x)\in \{ 0,r\} $
and
$\dualrank (x)\in \{ 0,d\} $.

In case $d=r=c$, suppose, for a contradiction, that $x$
is nonextremal.
Then
$\rank (x)\in \{ 0,r\} =\{ 0,c\} $ implies
$\rank (x)=c$
and
$\dualrank (x)\in \{ 0,d\} =\{ 0,c\} $ implies
$\dualrank (x)=c$.
However, this implies
$c=\rank (x)+\dualrank (x)=2c$, a contradiction.
Thus, in case $d=r=c$, we have that $x$ is maximal or minimal, and
we are left with the case that one of $d$ and $r$ is not $c$
and hence strictly smaller than $c$.

If $c>d\geq \dualrank (x)$, then $x$ is not minimal.
Hence the rank of $x$ must be $r$.
Now $x$ is nonextremal iff $r<c$.

If $c>r\geq \rank (x)$, then $x$ is not maximal.
Hence the dual rank of $x$ must be $d$.
Now,
the rank of $x$ is $c-d$, and
$x$ is nonextremal iff $d<c$.
\qed

\vspace{.1in}

{\bf Proof of Theorem \ref{getnonexcndeck}.}
By Proposition \ref{getranksofNTMAcards} and Theorem \ref{nonextcardthmnew},
Theorem \ref{getnonexcndeck} is proved for indecomposable ordered sets and
for ordered sets with $3$ or more elements in nontrivial order-autonomous subsets.
This leaves the case of
connected coconnected ordered sets $P$
with
exactly two elements in nontrivial order-autonomous subsets.
Let $A=\{ a_1 , a_2 \} $ be the unique
nontrivial order-autonomous subset of $P$
and note that the two NTMA-cards $P\setminus \{ a_1 \} $
and $P\setminus \{ a_2 \} $
are isomorphic.

By Theorem \ref{nonextcardthmnew}, all non-NTMA cards
have been identified as extremal or nonextremal.
For every $r>0$ let $n_r $ be the number of
nonextremal elements of rank $r$ in $P$ and let
$i_r $ be the number of nonextremal non-NTMA-cards
for which the rank of the removed element has been identified to be $r$.
Then $n_r -i_r \in \{ 0,1,2\} $.
In case
$n_r -i_r =0$, the deck
${\cal N} _P ^r $ has been reconstructed
as the multiset of nonextremal non-NTMA-cards
whose removed element has rank $r$.
In case
$n_r -i_r =1$, the deck
${\cal N} _P ^r $
consists of the
nonextremal non-NTMA-cards of rank $r$ plus one copy of
$P\setminus \{ a_1 \} $.
In case
$n_r -i_r =2$, the deck
${\cal N} _P ^r $
consists of the
nonextremal non-NTMA-cards of rank $r$ plus two copies of
$P\setminus \{ a_1 \} $.
\qed

\begin{remark}

{\rm
Note that the above also allows us to
reconstruct the
decks
${\cal NTMA} _P ^r :=\langle P\setminus \{ x\} :
x\in P, {\rm rank} (x)=r, P\setminus \{ x\} {\rm \ is \ an \ NTMA-card} \}
\rangle $
of unlabelled
NTMA-cards obtained by removing an element of rank $r\geq 0$
in the case that there are exactly
$2$ elements in nontrivial order-autonomous
subsets, which was not considered in
Proposition \ref{getranksofNTMAcards}:
The decks of nonextremal
cards (NTMA or not)
of rank $r>0$
have been identified in the proof of Theorem \ref{getnonexcndeck}
above.
In the identifiable case that
the set $A$ is a chain $\{ a_1 <a_2 \} $,
at most one card $P\setminus \{ a_i \} $ has not been
added to any of these decks.
If so,
it is
identifiable via
part \ref{nhoodfacts3} of Lemma \ref{nhoodfacts} and its dual
whether
this card
$P\setminus \{ a_i \} $ is a maximal or a minimal card.
If $A$ is an antichain with $2$ elements, then
$P\setminus \{ a_i \} $ is minimal (maximal) iff it has fewer
minimal (maximal) elements than $P$.
If an $a_i $ is minimal, its rank is $0$.
If an $a_j $ is maximal,
its rank is the height of the unique missing ideal from
the ideal deck of $P$
(see part \ref{nhoodfacts2} of Lemma \ref{nhoodfacts})
that is not in the ideal deck of $P\setminus \{ a_j \} $.

In particular, the above allows us to identify the
deck of minimal NTMA-cards and
the deck of maximal NTMA-cards.
Similar to the
proof of Theorem \ref{getnonexcndeck},
the
(still unfinished)
reconstruction of the
minimal and maximal decks can thus focus on non-NTMA cards.

}

\end{remark}

\begin{cor}
\label{getsomemincards}

Let $P$ be a connected ordered set with $\geq 4$ elements
and let $P\setminus \{ x\} $ be an extremal
non-NTMA card.
Unless
$P\setminus \min (P) $ has a component
$C$
that has a
minmax pair $(l,h)$ of pseudo-similar points
or that is a singleton $\{ l\} $, and
$x$ is the unique lower cover
of the element $l$,
we can identify whether $P\setminus \{ x\} $ is a minimal card.
Moreover, if $C$ is a chain, then $C$ is order-autonomous in $P$.
In this case, we call $C$ a
{\bf ranging chain}.

\end{cor}

{\bf Proof.}
In the proof of Theorem \ref{nonextcardthmnew}
above, if
one of $r$ or $d$ is not equal to the length of the longest
chain that contains $x$, then $x$ can be identified
as non-maximal or non-minimal and then subsequently as
minimal or maximal or nonextremal.
This leaves the case that both $r$ and $d$ are
equal to the length of the longest
chain that contains $x$.

If $x$ is minimal and $P\setminus \{ x\} $ has as many
minimal elements as $P$, then, in the proof of
Lemma \ref{getranknonex1} above,
$U_x :=\left( P\setminus \{ x\} \right) \setminus
\min \left( P\setminus \{ x\} \right) $
is a minimal card of $Q:=P\setminus \min (P)$.
When the $r$ from the
proof of
Lemma \ref{getranknonex1} above is the length of a longest chain that contains
$x$, then $U_x $ is also a maximal card of $Q$.
Consequently,
$Q$ has
a component
$C$ that has a
minmax pair $(l,h)$ of pseudo-similar points
or that is a singleton $\{ l\} $, and
$x$ is the unique lower cover
of the element $l$.

Moreover, if $C$ is a chain
$\{ c_1 <\cdots <c_n \} $,
then, because
the proof of
Lemma \ref{getranknonex1} leads to
$r=n$,
for any element $c_i <c_n $ of $C$, we must have that
$|\downarrow c_i |=|\downarrow c_{i+1} |-1$, which means that
$C$ is order-autonomous in $P$.
\qed

\section{Theorem \ref{minmaxpsrecon}:
A Consequence of Nonreconstructibility}
\label{conswhenconn}

The simplest, though naive, idea for reconstruction
of ordered sets with a minmax pair of pseudo-similar points
is to try to identify a minimal card that is isomorphic to a
maximal card and then use
Theorem \ref{keyiso}.
Lemma \ref{identnonbotofSUMRC} identifies
some maximal (and, dually, minimal) cards
in the extremal deck
identified by Theorem \ref{getnonexcndeck}.
However, if $(l,h)$ is a minmax pair of pseudo-similar points, then
$h$ (as well as,
in certain examples,
more maximal elements) does not
satisfy the hypothesis of
Lemma \ref{identnonbotofSUMRC}; and further
attempts, say, with Corollary \ref{getsomemincards}
and its dual,
also do not allow for the
identification of a minimal card that is isomorphic to a maximal card.

Therefore, to prove the reconstructibility of
ordered
sets with a minmax pair $(l,h)$ of pseudo-similar points,
we will assume that there
is a set $Q$
that has
the same deck as a set
$N$ with a minmax pair $(l,h)$ of pseudo-similar points
and that is not isomorphic to $N$.
From this assumption, Lemma \ref{dnewisolem} below
derives the existence of a
certain isomorphism between
maximal cards of
an ordered set
with a minmax pair
of pseudo-similar points.
Remark \ref{whytechnical} below indicates why
Lemma \ref{dnewisolem}
must refer to the action of the isomorphism on specific points
from Definition \ref{defineAP} below.
In Section \ref{isobetwcards}, we will then prove that
an isomorphism as claimed in Lemma \ref{dnewisolem} cannot exist.
Recall that a {\bf filter} is a subset $\uparrow x$
and that, by
the dual of part \ref{nhoodfacts2} of Lemma \ref{nhoodfacts},
the filter deck is reconstructible.

\begin{define}

Let $P$ be an ordered set
and
let $x\in P$ be maximal.
Then $x$ is called
{\bf filter shifting} iff the
filter deck
of
$C\setminus \{ x\} $
is obtained from the filter
deck of $P$ by removal of a
single filter.

\end{define}

\begin{lem}
\label{identnonbotofSUMRC}

Let $P$ be an ordered set
and let $C\setminus \{ x\} $
be an extremal card.
If $x$ is maximal and not filter shifting, then
$C$ is identifiable as a maximal card.

\end{lem}

{\bf Proof.}
The filter deck of
any minimal card
is obtained from the filter deck of $P$ by removal of a
single filter.
Because
$x$ is not filter shifting,
the filter deck of
$C$ is not
obtained from the filter deck of $P$ by removal of a
single filter.
Consequently, we can identify $C$
as a maximal card.
\qed

\begin{define}

Let $P$ and $S_1 , \ldots , S_k $
be ordered sets.
The point $c\in P$ is called a {\bf cutpoint}
iff $P\setminus \{ c\} $ has more components than $P$.
We say that removal of $c$ {\bf induces components} isomorphic to
$S_1 , \ldots , S_k $ iff
$c$ is an element of a component $K$ of $P$ such that
$K\setminus \{ c\} $ has at least $k$ components, and,
among the components of
$K\setminus \{ c\} $, there are pairwise distinct components
$C_1 , \ldots , C_k $
such that $C_i $ is isomorphic to $S_i $.

\end{define}

Let $(l,h)$ be a minmax pair of pseudo-similar points in $P$.
Because $K_h \supseteq \{ h\} \cup \bigcup _{j=1} ^r R_j $, clearly,
for all $j$, we have
$|K_h |>|R_j |$. Thus, we will refer to
$R_1 , \ldots , R_r $ as the {\bf small components}
of $P\setminus \{ h\} $
and to
$L_1 , \ldots , L_z $ as the {\bf small components}
of $P\setminus \{ l\} $.

\begin{lem}
\label{pdisclem}

Let $P$ be a connected
ordered set with a minmax pair $(l,h)$ of pseudo-similar points
such that $l$ and $h$ are cutpoints.
Then $P$ has a card that does not have a maximal element whose
removal induces components isomorphic to
the small components.
Dually, $P$ has a card that does not have a minimal element whose
removal induces components isomorphic to
the small components.

\end{lem}

{\bf Proof.}
Suppose for a contradiction that $P$ is a smallest
connected
ordered set with a minmax pair $(l,h)$ of pseudo-similar points
such that $l$ and
$h$ are cutpoints and
such that every card has a maximal element whose removal induces components
isomorphic to
the small components.
By part \ref{stuffaboutlh} of Lemma \ref{propertysequence},
no two components of $P\setminus \{ h\} $ have the same size.
Hence
we can assume that $|R_1 |>|R_j |$ for all $j\in \{ 2, \ldots , r\} $.
Let $y\in R_1 $ such
that $y$ is not a cutpoint of $P$
(any point in $R_1 $ at maximum distance from $h$ is not a cutpoint)
and let $K:=P\setminus \{ y\} $.
By assumption, $K
$ has a maximal element $m$
such that
$K\setminus \{ m\}$ has components
$A_1 , \ldots , A_r $ isomorphic to
$L_1 , \ldots , L_r $.
Because the components of $K\setminus \{ h\} $ are
$K_l $, $R_1 \setminus \{ y\} $ and $R_2 , \ldots , R_r $,
none of which is isomorphic to $R_1 $, we conclude that
$m\not= h$.

Now suppose for a contradiction that
$m\in L_j $ for some $j\in \{ 1, \ldots , r\} $.
Removal of $m$ induces components
$A_1 , \ldots , A_r $ isomorphic to $L_1 , \ldots , L_r $.
The component $A_j $ would not be contained in $L_j $, so $A_j $ would contain
$l$ and hence it would contain
$\{ l\} \cup C\cup \{ h\} \cup R_1 \setminus \{ y\} $.
However the latter set has
(even in case $j=1$ and $C=\emptyset $) at least
$|R_1 |+1> |A_j |$ elements, contradiction.
Thus $m\not\in L$ and similarly, $m\not\in R$.
Hence
$m\in C$.

Now
$h\in A_i $ would imply
that there is a lower bound $z$ of $h$
such that
$z\in C$ and
$\{ h,z\} \cup
R_1 \setminus \{ y\} \subseteq A_i $,
which is not possible,
because the former set has more elements than the latter set.
Hence $h\not\in A_i $.
Now,
$A_i \cap R_j \not= \emptyset $
would imply that there is
a path
from an element of $R_j $ through $A_i $ to $m\in C$.
Because this path would not contain $h$, this is not possible,
and hence, for all $i,j$, we have
$A_i \cap (R_j \cup \{ h\} )= \emptyset $.
Similarly, we prove
$l\not\in A_i $
($l\in A_i $ would imply $\{ l\} \cup L_i \subseteq A_i $, which is not
possible), and then
$A_i \cap (L_j \cup \{ l\} )=\emptyset $ for all $i,j$.

Thus $m$ is a cutpoint in $C\subseteq K_h \subseteq P\setminus \{ l\} $
whose removal induces components
$A_1 , \ldots , A_r \subseteq C$ isomorphic to
$L_1 , \ldots , L_r $, and $K_h $
has a minmax pair $(\Phi (l),h)$ of pseudo-similar points.
Because $\{ m\} \cup \bigcup _{j=1} ^z A_j \subseteq C$,
$C$ cannot be a small component of $K_h \setminus \{ h\} =
C\cup \bigcup _{j=1} ^z R_j $.
Hence $R_1 , \ldots , R_z $ are the
small components of $K_h \setminus \{ h\} $.

However, now
$K_h $ has fewer elements than $P$,
and every card of $K_h $ has a maximal element
($m$ or $h$, respectively)
whose removal induces components
that are isomorphic to the small components $R_1 , \ldots , R_z $.
This is a contradiction to the choice of $P$.

The dual claim follows by self-duality,
see part \ref{dualautlh} of Lemma \ref{propertysequence}.
\qed

\vspace{.1in}

For the next definition, recall that, for all
$x$ in an ordered set with a minmax pair $(l,h)$ of pseudo-similar
points, we have $|\downarrow x|\leq |\downarrow \Phi (x)|$.

\begin{define}
\label{defineAP}

Let $P$ be a connected ordered set with a minmax pair $(l,h)$ of pseudo-similar
points. We define $A_P \subseteq P$ to be the set of all maximal elements $x$ of $P$
such that
$|\downarrow x|=|\downarrow h|$.
Moreover, $n:=|P|$
and
$v$ is chosen so that
$A_P =\left\{ \Phi ^{v } (l), \Phi ^{v +1} (l), \ldots , \Phi ^{n-1} (l) \right\} $.

\end{define}

\begin{lem}
\label{abovel}

Let $P$ be a connected ordered set with a minmax pair $(l,h)$ of pseudo-similar
points
and let
$A_P =\left\{ \Phi ^{v } (l), \Phi ^{v +1} (l), \ldots , \Phi ^{n-1} (l) \right\} $.
Then
$\left| \uparrow \Phi ^j (l)\right| = |\uparrow l|$
iff $j\in \{ 0, \ldots , n-1-v \} $.

Moreover,
for all $j,k\in \{ 0, \ldots , n-1-v \} $, we have that
$\Phi ^j (l)\leq \Phi ^{v +j} (l)$
and that
$\Phi ^j (l)\leq \Phi ^{v +k} (l)$
implies $k\leq j$.

\end{lem}

{\bf Proof.}
Because $|A_P |=n-v$, by part \ref{dualautlh} of Lemma \ref{propertysequence},
$P$ has exactly
$n-v$ elements with $|\uparrow l|$ upper bounds.
Because, for all $x\in P$, we have
$|\uparrow x|\geq |\uparrow \Phi (x)|$,
these elements must be
$l, \Phi (l), \ldots , \Phi ^{n-1-v} (l)$.

{\em Claim: $l<\Phi ^v (l)$ and $\Phi ^v (l)$ is the only element of $A_P $ above $l$.}
Because $\left| \downarrow \Phi ^{v  } (l)\right| =|\downarrow h|
>
\left| \downarrow \Phi ^{v  -1} (l)\right|
=
\left| \downarrow \Phi ^{-1} \left( \Phi ^{v} (l)\right) \right|
$, we must have that
$\Phi ^{v  } (l)>l$.
Now suppose, for a contradiction, that there is an $m\in \{ v +1, \ldots , n-1\} $
such that
$\Phi ^{m} (l)>l$.
Then
$
|\downarrow h|
=
\left| \downarrow \Phi ^{m} (l)\right|
>
\left| \downarrow \Phi ^{-1} \left( \Phi ^{m} (l)\right) \right|
=
\left| \downarrow \Phi ^{m-1} (l)\right|
\geq \left| \downarrow \Phi ^{v  } (l)\right|
=|\downarrow h|$, a contradiction.

Now let $j\in \{ 0, \ldots , n-1-v \} $.
Application of $\Phi ^j $
to
$l<\Phi ^{v  } (l)$
shows that
$\Phi ^j (l)\leq \Phi ^{v +j} (l)$.
Moreover, suppose that there was a $k\in \{ j+1, \ldots , n-1-v \} $ such that
$\Phi ^j (l)\leq \Phi ^{v +k} (l)$. Then
we would have
$l
=\Phi ^{-j} \left( \Phi ^j (l)\right)
\leq
\Phi ^{-j} \left( \Phi ^{v +k} (l)\right)
= \Phi ^{v +k-j} (l)$, contradicting the {\em Claim}.
Hence $\Phi ^j (l)\leq \Phi ^{v +k} (l)$ implies $k\leq j$.
\qed

\begin{define}

Let $P$ be a connected ordered set with a minmax pair $(l,h)$ of pseudo-similar
points.
For each $p=\Phi ^{v+j} (l) \in A_P$,
we set
$d_p =\Phi ^j (l)$.

\end{define}

\begin{lem}
\label{dnewisolem}

Let $N$ be a connected ordered set
with a minmax pair $(l_N , h_N )$ of pseudo-similar points.
Suppose there is a
connected ordered set $Q$
that
is not isomorphic to $N$, but that
has the same deck as $N$.
Then there is a
connected ordered set $P
$
with a minmax pair
$(l,h)$ of pseudo-similar
points
and not necessarily distinct
maximal elements
$a,b\in A_P $
such that
there is an isomorphism $\Psi :P \setminus \{ a\} \to P \setminus \{ b\} $
that satisfies
$\Psi (d_a )\not= d_b $, and such that
each of $a$ and $b$ has
exactly one
lower bound with
$\left| \uparrow l \right| $ upper bounds in $P$.

\end{lem}

{\bf Proof.}
Because $Q$ and $N$ have the same deck, they have
equal ideal size sequences.
Because two ordered sets with
minmax pairs of pseudo-similar points and
the same ideal size sequences
must by Theorem \ref{keyiso}
be isomorphic, $Q$ does not have any
minmax pairs of pseudo-similar points.
By
Theorem \ref{getnonexcndeck}
and duality, we can
therefore assume without loss of generality that
$Q$ has maximal elements $a$ and $x$ such that the cards
$Q\setminus \{ a\} $, $Q\setminus \{ x\} $,
$N\setminus \{ l_N \} $ and $N\setminus \{ h_N \} $
are isomorphic
and such that
at least one
of
$Q\setminus \{ a\} $ and $Q\setminus \{ x\} $,
is not identifiable as maximal.
Consequently, because they are isomorphic,
neither of
$Q\setminus \{ a\} $ and $Q\setminus \{ x\} $,
is identifiable as maximal.

Suppose, for a contradiction that the card
$N\setminus \{ h_N \} $ is disconnected with
components
$K_{h_N} , R_1 , \ldots , R_z $, $z\geq 1$, and
$\left| K_{h_N } \right| >|R_j |$ for all $j$.
Then
$Q\setminus \{ a\} $
has components
$K^a , R_1 ^a, \ldots , R_z ^a $,
and
$Q\setminus \{ x\} $
has components
$K^x , R_1 ^x, \ldots , R_z ^x $,
with notation such that each component is isomorphic to the
corresponding component of
$N\setminus \{ h_N \} $.
Now $x\in R_i ^a $ would imply that
$Q\setminus \{ x\} $ has a component that contains
$K^a \cup \{ a\} $, which has more elements
than any component of
$Q\setminus \{ x\} $.
Thus $x\in K^a $.
Similarly $a\in K^x $
and
then
$R_1 ^x , \ldots , R_z ^x \subseteq K^a $
and
$R_1 ^a , \ldots , R_z ^a \subseteq K^x $.
In particular, we obtain $Q\subseteq K^a \cup K^x $.

Now let $q\in Q$.
If $q\in K^x $, then
removal of $x$ induces components
$R_1 ^x , \ldots , R_z ^x $
isomorphic to
$R_1 , \ldots , R_z $.
If $q\in K^a $, then
removal of $a$ induces components
$R_1 ^a , \ldots , R_z ^a $
isomorphic to
$R_1 , \ldots , R_z $.
Therefore, every card of $Q$
has a maximal cutpoint
$c$
such that
removal of $c$ induces components isomorphic to
$R_1 , \ldots , R_z $.
However, by Lemma \ref{pdisclem}, $N$ has
a card that
has no maximal element
whose removal induces components isomorphic to
$R_1 , \ldots , R_z $, contradicting that
$Q$ and $N$ have isomorphic decks.
Therefore the cards
$Q\setminus \{ a\} $, $Q\setminus \{ x\} $,
$N\setminus \{ l_N \} $ and $N\setminus \{ h_N \} $
are connected.

Because
$|\downarrow a|=|\downarrow h_N |$
and
(clearly) $x\not<a $,
we have that $a$ is maximal in $P_a :=Q\setminus \{ x\} $
with $|\downarrow
_{P_a }
a|=|\downarrow h_N |$ lower bounds.
Note that,
because it is isomorphic to $N\setminus \{ l_N \} $,
$P_a $ is connected and
has a minmax pair of pseudo-similar points.
Because $Q$ has the same ideal size sequence as $N$, no element of $Q$ has more than
$|\downarrow h_N |$ lower bounds, every element of
$A
_{P_a } $ has
$|\downarrow h_N |$ lower bounds
and, for every $z\in A
_{P_a }
$
the element $d_z $ has
$|\uparrow l_N |$ upper bounds in
$P_a $.

Consider the element $d_a \in P_a $.
By the above, $d_a
$ has $|\uparrow l_N |$
upper bounds in $P_a $.
If we had $x>d_a $, then, in $Q$, the element $d_a $ would have
$|\uparrow l_N |+1$ upper bounds.
However,
because $Q$ has the same filter size sequence as $N$,
no element of $Q$ has more than
$|\uparrow l_N |$ upper bounds.
Thus $x\not> d_a $.

Let
$P_x :=Q\setminus \{ a\} $,
let
$\Gamma :P_a \to P_x $
be an isomorphism,
and
let $b:=\Gamma ^{-1} (x)
$.
Because
$P_a \setminus \{ a\}
=Q\setminus \{ a,x\}
=P_x \setminus \{ x\} $,
for the ordered set $P_a $,
the maximal card
$P_a \setminus \{ a\} $
is isomorphic to
the maximal card
$
\Gamma ^{-1} [P_x \setminus \{ x\} ]
=
P_a \setminus \left\{ \Gamma ^{-1} (x)\right\}
=
P_a \setminus \left\{ b\right\}
$.
Let
$\Psi :P_a \setminus \{ a\} \to P_a\setminus \{ b\} $
be an isomorphism.
In case $\Psi (d_a )\not= d_b $, the isomorphism
$\Psi :P_a \setminus \{ a\} \to P_a \setminus \{ b\} $
satisfies the required inequality.

This leaves the case in which
$\Psi (d_a )=d_b $.
In this case,
$\widetilde{\Psi }:=\Gamma \circ \Psi $
is an isomorphism from $P_a \setminus \{ a\} $ to
$\Gamma [\Psi [P_a \setminus \{ a\} ] ]
=
\Gamma [P_a \setminus \{ b\}  ]
=
\Gamma \left[ P_a \setminus \left\{ \Gamma ^{-1} (x)\right\} \right]
=
P_x \setminus \{ x\}
=
P_a \setminus \{ a\}
$.
Because, by
part \ref{lharerigid} of Lemma \ref{propertysequence},
$P_a $ and $P_x $ are rigid,
we have
$
\widetilde{\Psi } (d_a)=
\Gamma (\Psi (d_a ))=\Gamma (d_b )=\Gamma (d_{\Gamma ^{-1} (x)} )=d_x $,
which, because $x\not> d_a $ is not
equal to $d_a $.
Therefore, in case $\Psi (d_a )=d_b $, with
$\widetilde{b}:=a$, the isomorphism
$\widetilde{\Psi } :P_a \setminus \{ a\}
\to P_a \setminus \left\{ \widetilde{b} \right\} $
satisfies the required inequality.

If necessary, rename $\widetilde{\Psi } $ to $\Psi $ and
$\widetilde{b } $ to $b$.

Finally,
because
$Q\setminus \{ a\} $ and $N\setminus \{ h_N \} $
are isomorphic, we have
$|\downarrow _{P_a } a|=|\downarrow _Q a|=|\downarrow h_N |$,
and, because no element of $Q$ can have more lower bounds,
$a\in A_{P_a } $.
Isomorphism of
$P_a \setminus \{ a\} $ and $P_a \setminus \{ b\} $
and maximality of $b$
imply $b\in A_{P_a} $.

Because $Q\setminus \{ a\} $ is
not identifiable as a maximal card, by Lemma \ref{identnonbotofSUMRC},
$a$ must be filter shifting, which means
$a$ has at most one
lower bound with
$\left| \uparrow l_N \right| $ upper bounds
in $Q$.
Consequently,
$a$ has at most
one lower bound with
$\left| \uparrow l_N \right| $ upper bounds in $P_a  \subseteq Q$.
Because $a$ has
$\left| \downarrow h_N \right| =\left| \uparrow l_N \right| $
lower bounds in $P_a $,
by Lemma \ref{abovel},
it has
at least (and hence exactly) one
lower bound with
$\left| \uparrow l_N \right| $ upper bounds in $P_a $.

Similarly,
$x$ has
exactly one
lower bound with
$\left| \uparrow l_N \right| $ upper bounds in $P_x $
and hence
$\Gamma ^{-1} (x)$ has
exactly one
lower bound with
$\left| \uparrow l_N \right| $ upper bounds in $P_a $.

Because $b$ is either $\Gamma ^{-1} (x)$
or $a$,
it has
exactly one
lower bound with
$\left| \uparrow l_N \right| $ upper bounds in $P_a $, too.
\qed

\begin{figure}

\centerline{
\unitlength 1mm 
\linethickness{0.4pt}
\ifx\plotpoint\undefined\newsavebox{\plotpoint}\fi 
\begin{picture}(115.5,27)(0,0)
\put(5,5){\line(0,1){20}}
\put(50,5){\line(0,1){20}}
\put(95,5){\line(0,1){20}}
\put(5,25){\line(1,-1){20}}
\put(50,25){\line(1,-1){20}}
\put(95,25){\line(1,-1){20}}
\put(25,5){\line(0,1){20}}
\put(70,5){\line(0,1){20}}
\put(115,5){\line(0,1){20}}
\put(5,5){\circle*{2}}
\put(50,5){\circle*{2}}
\put(95,5){\circle*{2}}
\put(25,5){\circle*{2}}
\put(70,5){\circle*{2}}
\put(115,5){\circle*{2}}
\put(5,15){\circle*{2}}
\put(50,15){\circle*{2}}
\put(95,15){\circle*{2}}
\put(25,15){\circle*{2}}
\put(70,15){\circle*{2}}
\put(115,15){\circle*{2}}
\put(5,25){\circle*{2}}
\put(50,25){\circle*{2}}
\put(95,25){\circle*{2}}
\put(25,25){\circle*{2}}
\put(70,25){\circle*{2}}
\put(115,25){\circle*{2}}
\put(25,5){\line(5,4){25}}
\put(70,5){\line(5,4){25}}
\put(50,27){\makebox(0,0)[cb]{\footnotesize $a$}}
\put(25,3){\makebox(0,0)[ct]{\footnotesize $l=d_a $}}
\put(95,27){\makebox(0,0)[cb]{\footnotesize $h$}}
\end{picture}
}

\caption{
A set with a minmax pair of pseudo-similar points
and a non-rigid
maximal card $P\setminus \{ a\} $
obtained by removing a
maximal element $a$ with $|\downarrow h|$ lower bounds.
(The component
of $P\setminus \{ a\} $ that contains $h$ is not rigid.)
}
\label{whyneedca}

\end{figure}
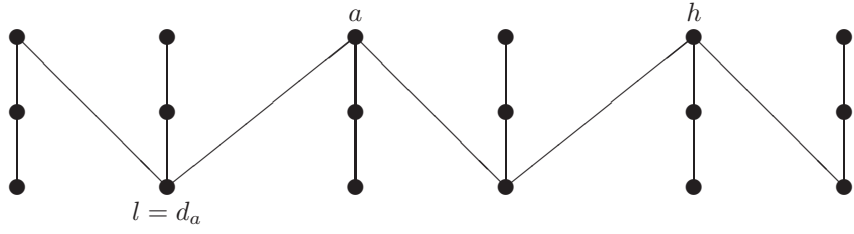

\begin{remark}
\label{whytechnical}

{\rm
The technical
condition $\Psi (d_a )\not= d_b $
assures that, in case $a=b$,
$\Psi $ is a nontrivial automorphism.
Figure \ref{whyneedca} shows
that
maximal cards can
have
nontrivial automorphisms.
Therefore, in the proof that there are no sets
$P$ as
specified in Lemma \ref{dnewisolem},
we must explicitly consider the
element $d_a $.
}

\end{remark}

\section{Theorem \ref{minmaxpsrecon}:
Properties of Isomorphisms Between Cards}
\label{isobetwcards}

This section investigates further
properties that
an isomorphism
$\Psi $ as in Lemma \ref{dnewisolem}
must have.
Recall that $R:=\bigcup _{j=1} ^z R_j $
denoted the union of the small components.

\begin{lem}
\label{sameorderofremoval}

Let $P$ be a connected ordered set with
a minmax pair $(l,h)$ of pseudo-similar
points
and let $\Psi :K_l \setminus \left\{ \Phi ^{-1} (h)\right\} \to K_l \setminus \{ l\} $
be the
(by
parts \ref{stuffaboutlh} and \ref{lharerigid}
of Lemma \ref{propertysequence})
unique isomorphism between these two sets.
For every
$x\in K_l \setminus \left\{ \Phi ^{-1} (h)\right\} $, define
$$
\widehat{\Phi } (x)
:=
\cases{
\Phi (x) & if $x\not\in \Phi ^{-1} [R]$, \cr
\Phi ^2 (x) & if $x\in \Phi ^{-1} [R]$. \cr
}
$$
Then $\Psi = \widehat{\Phi } $.

\end{lem}

{\bf Proof.}
Let
$x\in K_l \setminus \left\{ \Phi ^{-1} (h)\right\} $.
Then $\Phi (x)\not\in \{ l,h\} $.
If
$x\not\in \Phi ^{-1} [R]$, then
$\widehat{\Phi } (x) = \Phi (x)\not \in R\cup \{ h\} $ and
hence $\widehat{\Phi } (x)\in K_l \setminus \{ l\} $.
If
$x\in \Phi ^{-1} [R]$, then
$\widehat{\Phi } (x) = \Phi ^2 (x)\in \Phi [R]=L\subseteq K_l \setminus \{ l\} $.
Hence
$\widehat{\Phi } \left[
K_l \setminus \left\{ \Phi ^{-1} (h)\right\}
\right]
\subseteq
K_l \setminus \{ l\} $.

Next, we prove
that
$x\leq y$ iff $\widehat{ \Phi } (x) \leq \widehat{\Phi } (y)$.
To this end, let $x,y\in K_l \setminus \left\{ \Phi ^{-1} (h)\right\} $.
If either $x,y\in \Phi ^{-1} [R]$
or
$x,y\not\in \Phi ^{-1} [R]$,
then, because $\Phi $ and $\Phi ^2 $ are both isomorphisms,
we have $x\leq y$ iff $\widehat{ \Phi } (x) \leq \widehat{\Phi } (y)$.
Now consider the case that
$x\in \Phi ^{-1} [R]$ and $y\not\in \Phi ^{-1} [R]$.
First suppose, for a contradiction, that $x$ is comparable to $y$.
Then $\Phi (x)\in R$ is comparable to $\Phi (y)\not\in R$,
which is only possible if $\Phi (y)=h$.
However, this contradicts
$y\in K_l \setminus \left\{ \Phi ^{-1} (h)\right\} $.
Hence $x $ is not comparable to
$y$.
Now suppose,
for a contradiction, that
$\widehat{\Phi } (x) $ is comparable to $\widehat{\Phi } (y)$.
Then
$\Phi ^2 (x) =\widehat{\Phi } (x) $ is comparable to $\Phi (y)=\widehat{\Phi } (y) $,
which means that
$\Phi (x)\in R$ is comparable to $y\in K_l \setminus \left\{ \Phi ^{-1} (h)\right\} $,
which is not possible, because removal of $h$ separates $K_l $ and $R$.
Hence
$\widehat{\Phi } (x) $ is not comparable to $\widehat{\Phi } (y)$.
Consequently, in case
$x\in \Phi ^{-1} [R]$ and $y\not\in \Phi ^{-1} [R]$,
the equivalence
$x\leq y$ iff $\widehat{ \Phi } (x) \leq \widehat{\Phi } (y)$
is vacuously true.
The
case
$x\not\in \Phi ^{-1} [R]$ and $y\in \Phi ^{-1} [R]$ is handled similarly.

In particular, $x\leq y$ iff $\widehat{ \Phi } (x) \leq \widehat{\Phi } (y)$
establishes that $\widehat{\Phi } $ is
injective.
Because
domain and range have the same number of elements
and $\widehat{\Phi } $ is
injective,
$\widehat{\Phi } $ is
bijective and hence
it is an isomorphism from $K_l \setminus \left\{ \Phi ^{-1} (h)\right\} $
to $K_l \setminus \{ l\} $.

By parts \ref{stuffaboutlh} and \ref{lharerigid}
of Lemma \ref{propertysequence},
$K_l \setminus \left\{ \Phi ^{-1} (h)\right\} $
and $K_l \setminus \{ l\} $
are rigid, which
implies that $\widehat{\Phi } = \Psi $.
\qed

\begin{lem}
\label{findlnew}

Let $P$ be a connected ordered set
with a minmax pair $(l,h)$ of pseudo-similar
points, let
$H_P :=\left\{ \Phi  ^{v+w} (l) , \ldots , \Phi ^{n-1} (l) \right\} \subseteq A_P $,
and
let $B$ be the component of $P\setminus H_P $ that contains $l$.
Then
$B$ has a minmax pair $(l,x)$ of pseudo-similar points and
no component of $P\setminus H_P $ has more elements than $B$.
Moreover, there is an $m> 0$ such that the
components of $P\setminus H_P $ of size $|B|$ are
$B, \Phi [B], \ldots , \Phi ^{m-1} [B]$.

\end{lem}

{\bf Proof.}
We first prove
by induction on the size of $H_P $ that
$B$ has a minmax pair $(l,x)$ of pseudo-similar points.
To anchor the induction, first note that,
if $|H_P |=0$, then there is nothing to prove.
Moreover, if $|H_P |=1$, then $H_P =\{ h\} $,
$P\setminus H_P =P\setminus \{ h\} $,
and $B =K_l $.
Now, by part \ref{largcomparelh} of Lemma \ref{componiso},
$B$ has a minmax pair $(l,x)$ of pseudo-similar points.

Now let $P$ be a connected ordered
set with a minmax pair $(l,h)$ of pseudo-similar
points, let $|H_P |\geq 2$,
and assume that the conclusion holds for all
connected ordered sets
$\overline{P}$
with a minmax pair of pseudo-similar
points
and sets $H_{\overline{P} } $
such that $\left| H_{\overline{P} } \right| <|H_P |$.
Because $|H_P |\geq 2$,
we have
$\left\{ \Phi ^{-1} (h), h\right\} \subseteq H_P $.
By part \ref{largcomparelh} of Lemma \ref{componiso},
the set $K_l $ is
a connected ordered set with a minmax pair
$\left( l,\Phi ^{-1} (h)\right) $ of pseudo-similar
points.
Let $\widehat{H}_{K_l } :=H_P \cap K_l $.
Clearly,
$\left| \widehat{H}_{K_l } \right| < |H_P |$ and
$B $ is a component of
$K_l \setminus \widehat{H}_{K_l } $.

Let
$\Psi _{K_l } :K_l \setminus \{ \Phi ^{-1} (h)\} \to K_l \setminus \{ l\} $
be the unique isomorphism between these sets and let
$s$ be the smallest integer such that $\Psi ^s (l)\in \widehat{H}_{K_l }$.
By Lemma \ref{sameorderofremoval},
every $\Psi ^{s+k} (l)$ is of the form
$\Phi ^{v+w+k'} (l)$ and hence in $\widehat{H}_{K_l} $, too.
Thus
$\widehat{H}_{K_l } =\left\{ \Psi ^{s} (l), \ldots , \Psi ^{|K_l |-1} (l) \right\} $
and,
by induction hypothesis,
$B$ has a minmax pair $(l,x)$ of pseudo-similar points.
This completes the induction.

For every component $D$ of $P\setminus H_P $
that does not contain $l$,
we define $\Upsilon (D)$ to be the component of
$P\setminus H_P $ that contains
$\Phi ^{-1} [D]$.
Then
$\left| \Upsilon (D)\right|
\geq
\left| \Phi ^{-1} [D]\right|
=|D|$.
For every component $D$ of $P\setminus H_P $,
we define
$j_D :=\min \left\{ j:\Phi ^j (l)\in D\right\} $.
Let $D$
be a
component of $P\setminus H_P $
that does not contain $l$.
Then
there is a $k\in \{ 1, \ldots , j_D \} $
such that
$\Upsilon ^k (D)=B$.
In particular, we obtain
$|D|\leq |B|$.

To prove the final claim,
let $D$
be a
component of $P\setminus H_P $
that does not contain $l$ and such that
$|D|=|B|$.
Then
$
|B|\geq
\left| \Upsilon (D)\right|
\geq
\left| \Phi ^{-1} [D]\right|
=|D|=|B|$, which implies
$\left| \Upsilon (D)\right|
=
\left| \Phi ^{-1} [D]\right|
$ and hence
$\Upsilon (D)
=
\Phi ^{-1} [D]
$.
In particular, this implies that
$j_{\Upsilon (D)} =j_D -1$.

Let $E$ be the
component of $P\setminus H_P $
with $|E|=|B|$ such that
$j_E$ is as large as
possible.
By the above, for all $k\in \{ 0,\ldots , j_E \} $,
we have that
$\Upsilon ^k (E)$
is defined, has $|B|$ elements and
$j_{\Upsilon ^k (E)} = j_E -k$.
Because $j_E $
was chosen to be as large as possible,
the components of $P\setminus H_P $
of size $|B|$
are $\Upsilon ^k (E)=\Phi ^{j_E -k} [B]$,
with $k\in \{ 0,\ldots , j_E \} $, which proves the
final claim.
\qed

\begin{define}

Let $P$ be an ordered set with a minmax pair $(l,h)$ of
pseudo-similar points
and let
$A_P =\left\{ \Phi ^{v } (l), \Phi ^{v +1} (l), \ldots ,
\Phi ^{n-1} (l) \right\} $.
For every $k\in \{ 0,\ldots , n-1-v\} $,
we
call the components
of $P\setminus
\left\{ \Phi ^{v+k } (l), \Phi ^{v +k+1} (l), \ldots ,
\Phi ^{n-1} (l) \right\} $
that have as many elements
as the component that contains $l$
the {\bf large components}.

\end{define}

\begin{figure}

\centerline{
\unitlength 1mm 
\linethickness{0.4pt}
\ifx\plotpoint\undefined\newsavebox{\plotpoint}\fi 
\begin{picture}(134,22)(0,0)
\put(70,10){\line(0,1){5}}
\put(130,10){\line(0,1){5}}
\put(50,10){\line(0,1){5}}
\put(110,10){\line(0,1){5}}
\put(30,10){\line(0,1){5}}
\put(90,10){\line(0,1){5}}
\put(10,10){\line(0,1){5}}
\put(70,15){\line(-1,-1){5}}
\put(130,15){\line(-1,-1){5}}
\put(50,15){\line(-1,-1){5}}
\put(110,15){\line(-1,-1){5}}
\put(30,15){\line(-1,-1){5}}
\put(90,15){\line(-1,-1){5}}
\put(10,15){\line(-1,-1){5}}
\put(65,10){\line(0,1){5}}
\put(125,10){\line(0,1){5}}
\put(45,10){\line(0,1){5}}
\put(105,10){\line(0,1){5}}
\put(25,10){\line(0,1){5}}
\put(85,10){\line(0,1){5}}
\put(5,10){\line(0,1){5}}
\put(70,10){\circle*{1}}
\put(130,10){\circle*{1}}
\put(50,10){\circle*{1}}
\put(110,10){\circle*{1}}
\put(30,10){\circle*{1}}
\put(90,10){\circle*{1}}
\put(10,10){\circle*{1}}
\put(65,10.1){\circle*{1}}
\put(125,10.1){\circle*{1}}
\put(45,10.1){\circle*{1}}
\put(105,10.1){\circle*{1}}
\put(25,10.1){\circle*{1}}
\put(85,10.1){\circle*{1}}
\put(5,10.1){\circle*{1}}
\put(65,15){\circle*{1}}
\put(125,15){\circle*{1}}
\put(45,15){\circle*{1}}
\put(105,15){\circle*{1}}
\put(25,15){\circle*{1}}
\put(85,15){\circle*{1}}
\put(5,15){\circle*{1}}
\put(70,15){\circle*{1}}
\put(130,15){\circle*{1}}
\put(50,15){\circle*{1}}
\put(110,15){\circle*{1}}
\put(30,15){\circle*{1}}
\put(90,15){\circle*{1}}
\put(10,15){\circle*{1}}
\put(65,15){\line(1,1){5}}
\put(125,15){\line(1,1){5}}
\put(45,15){\line(1,1){5}}
\put(105,15){\line(1,1){5}}
\put(25,15){\line(1,1){5}}
\put(85,15){\line(1,1){5}}
\put(50,10){\line(-1,-1){5}}
\put(110,10){\line(-1,-1){5}}
\put(30,10){\line(-1,-1){5}}
\put(90,10){\line(-1,-1){5}}
\put(10,10){\line(-1,-1){5}}
\put(70,20){\circle*{1}}
\put(130,20){\circle*{1}}
\put(50,20){\circle*{1}}
\put(110,20){\circle*{1}}
\put(30,20){\circle*{1}}
\put(90,20){\circle*{1}}
\put(45,5){\circle*{1}}
\put(105,5){\circle*{1}}
\put(25,5){\circle*{1}}
\put(85,5){\circle*{1}}
\put(5,5){\circle*{1}}
\qbezier(70,20)(74,15)(70,10)
\qbezier(130,20)(134,15)(130,10)
\qbezier(50,20)(54,15)(50,10)
\qbezier(110,20)(114,15)(110,10)
\qbezier(30,20)(34,15)(30,10)
\qbezier(90,20)(94,15)(90,10)
\qbezier(45,5)(41,10)(45,15)
\qbezier(105,5)(101,10)(105,15)
\qbezier(25,5)(21,10)(25,15)
\qbezier(85,5)(81,10)(85,15)
\qbezier(5,5)(1,10)(5,15)
\put(45,3){\makebox(0,0)[t]{\footnotesize $l$}}
\put(90,22){\makebox(0,0)[cb]{\footnotesize $h$}}
\put(30,22){\makebox(0,0)[cb]{\footnotesize $\Phi ^{-1} (h)$}}
\put(110,22){\makebox(0,0)[cb]{\footnotesize $\Phi ^{-2} (h)$}}
\put(50,22){\makebox(0,0)[cb]{\footnotesize $\Phi ^{-3} (h)$}}
\put(5,5){\line(5,3){25}}
\put(25,5){\line(5,3){25}}
\put(85,5){\line(5,3){25}}
\put(45,5){\line(5,3){25}}
\put(105,5){\line(5,3){25}}
\put(70,10){\line(2,1){20}}
\end{picture}
}

\caption{
By Lemma \ref{findlnew},
no component of $P\setminus A_P $ is larger than the component $B $ that contains $l$.
The set above shows that there can, however,
be a component that has as many elements as
$B$.
}
\label{sevkll2}

\end{figure}
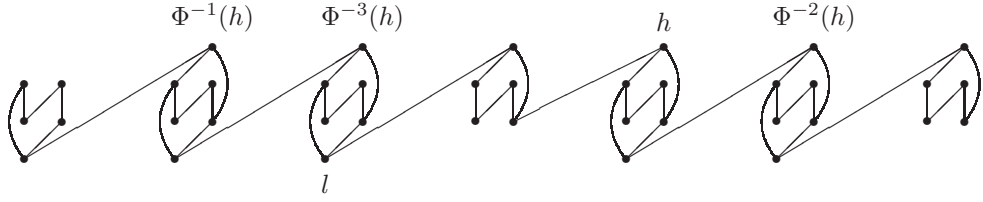

\begin{remark}

{\rm
The ordered set in Figure \ref{sevkll2} shows that
we can have more than one large component.
}

\end{remark}

Note that, by Lemma \ref{findlnew} and by
part \ref{lharerigid} of Lemma \ref{propertysequence},
the large components are rigid.

\begin{lem}
\label{partition}

Let $P$ be an ordered set with a minmax pair $(l,h)$ of
pseudo-similar points
and let
$A_P =\left\{ \Phi ^{v } (l), \Phi ^{v +1} (l), \ldots ,
\Phi ^{n-1} (l) \right\} $.
Then there are
$0=k_0 <k_1 <\cdots <k_X =n-v$ such that, for all
$x\in \{ 0,\ldots , X-1\} $,
the
large components
of $P\setminus
\left\{ \Phi ^{v +k_x } (l), \Phi ^{v +k_x +1} (l), \ldots ,
\Phi ^{n-1} (l) \right\} $
can be enumerated as
$B_{k_x } , \ldots , B_{k_{x+1}-1 } $
such that
the
following hold.
\begin{enumerate}
\item
\label{partition2}
For all $k\in \{ k_x , \ldots , k_{x+1} -2\} $, we have
$\Phi [B_k ]=B_{k+1} $.
\item
\label{partition3}
For all $k\in \{ k_x , \ldots , k_{x+1} -1\} $, we have
$\Phi ^k (l)\in B_k $.
\end{enumerate}

\end{lem}

{\bf Proof.}
We shall
assume that there are $M\geq 0$
and
$0=k_0 <k_1 <\cdots <k_M <n-v$ such that
parts
\ref{partition2} and \ref{partition3}
hold
for all
$x\in \{ 0,\ldots , M-1\} $
and such that
$l, \Phi ^{k_M } (l) \in B_{k_M } $.
Note that
the above is vacuously satisfied for $M=0$
and
$k_0 =0$ if we let
$B_0 $ be the component of
$P\setminus A_P $ that contains $l$.

By Lemma \ref{findlnew},
the large components of
$P\setminus
\left\{ \Phi ^{v +k_M } (l), \Phi ^{v +k_M +1} (l), \ldots ,
\Phi ^{n-1} (l) \right\} $
can be enumerated $B_{k_M } , \ldots , B_{k_M +m-1} $
such that, for all $k\in \{ k_M , \ldots , k_M +m-2\} $, we have
$\Phi [B_k ]=B_{k+1} $.
This establishes part \ref{partition2}
with $k_{M+1} :=k_M +m $.
Because $\Phi ^{k_M } (l) \in B_{k_M } $,
part \ref{partition3} is established, too.
In case $k_{M+1} =n-v$, the argument is complete
with $X:=M+1$.
Otherwise, let $B_{k_{M+1} } $ be the component of
$P\setminus
\left\{ \Phi ^{v +k_{M+1} } (l), \Phi ^{v +k_{M+1} +1} (l), \ldots ,
\Phi ^{n-1} (l) \right\} $
that contains $l$.
All that is left to do
is to
establish that
$
\Phi ^{k_{M+1} } (l) \in B_{k_{M+1} } $.

Clearly,
$N:=\Phi ^m [B_{k_M } ]$ is a component
of
$P\setminus \left\{ l,\Phi ^{v+{k_M }+1} (l), \ldots , \Phi ^{n-1} (l)\right\} $
with
$\left| B_{k_M } \right| $ elements that contains
$\Phi ^m \left( \Phi ^{k_M } (l)\right)
= \Phi ^{k_{M} +m} (l) = \Phi ^{k_{M+1} } (l) $.
If there was a $j\in \{ 0,\ldots , m-1\} $ such that
$N=\Phi ^{j } [B_{k_M } ]$,
then we would have $\Phi ^{m-j} [B_{k_M} ]=B_{k_M } $
and then, by rigidity, $\Phi ^{m-j} (l)=l$, which cannot be.
Hence,
$N$ is not equal to any $\Phi ^{k_M +j } [B_{k_M } ]$
with $j\in \{ 0,\ldots , m-1\} $.
Because these are the large components of
$P\setminus
\left\{ \Phi ^{v +k_M } (l), \Phi ^{v +k_M +1} (l), \ldots ,
\Phi ^{n-1} (l) \right\} $, we conclude that
$\Phi ^{v+k_M } (l)\in N$.
Because $\Phi ^{v+k_M } (l)>\Phi ^{k_M } (l)$
and $\Phi ^{k_M } (l)$ is, in
$B_{k_M } \subseteq P\setminus
\left\{ \Phi ^{v +k_M } (l), \Phi ^{v +k_M +1} (l), \ldots ,
\Phi ^{n-1} (l) \right\} $,
connected to $l$, we have that
$N\subseteq
B_{k_{M+1} } $
and in particular
$\Phi ^{k_{M+1} } (l) \in B_{k_{M+1} } $.
\qed

\begin{lem}
\label{partitionPsi}

Let $P$ be an ordered set with a minmax pair $(l,h)$ of
pseudo-similar points
and let
$A_P =\left\{ \Phi ^{v } (l), \Phi ^{v +1} (l), \ldots ,
\Phi ^{n-1} (l) \right\} $.
Let $a,b\in A_P $ be not necessarily distinct points
such that
each of $a$ and $b$ has
exactly one
lower bound with
$\left| \uparrow l \right| $ upper bounds in $P$
and
let
$\Psi : P\setminus \{ a\} \to P\setminus \{ b\} $ be an
isomorphism.
Then
$\Psi (d_a )=d_b $.

\end{lem}

{\bf Proof.}
Let $i_a , i_b $ such that
$a=\Phi ^{v+i_a } (l)$
and
$b=\Phi ^{v+i_b } (l)$, and
let $0=k_0 <k_1 <\cdots <k_X =n-v$ be as in
Lemma \ref{partition}.

{\em Claim.
If $k_{x+1} \leq \min \{ i_a , i_b \} $, then
$\Psi $ maps each of the sets
$\left\{ \Phi ^{k_x} (l) , \ldots , \Phi ^{{k_{x+1}}-1} (l)\right\} $
and
$\left\{ \Phi ^{v+k_x} (l) , \ldots , \Phi ^{v+{k_{x+1}}-1} (l)\right\} $
to itself.
}
First
note that the {\em Claim} is vacuously true for
$x=-1$.
Let $x\geq 0$ such that
$k_{x+1} \leq \min \{ i_a , i_b \} $
and assume the {\em Claim} has already
been established for all
$z\in \{ -1, \ldots , x-1\} $.
Then, because $\Psi $ maps $P\setminus A_P $ to itself
and
the {\em Claim} holds for all $z<x$,
$\Psi $ maps
$P\setminus \left\{ \Phi ^{v+k_x } (l), \ldots , \Phi ^{n-1} (l)\right\} $
to itself.
Let
the
large components
of $P\setminus
\left\{ \Phi ^{v +k_x } (l), \Phi ^{v +k_x +1} (l), \ldots ,
\Phi ^{n-1} (l) \right\} $
be enumerated
$B_{k_x } , \ldots , B_{k_{x+1}-1 } $
as in Lemma \ref{partition}.
Then $\Psi $ must be a permutation of the
rigid sets
$B_{k_x } , \ldots , B_{k_{x+1}-1 } $.

Let $i\in \{ 0,\ldots , k_{x+1} -1 -k_x \} $.
Then there is a
$j\in \{ 0,\ldots , k_{x+1} -1 -k_x \} $
such that
$\Psi \left[ B_{k_x +i } \right]
=B_{k_x +j } $.
Because the sets are rigid,
we obtain
$\Psi |_{B_{k_x +i } } =\Phi |_{B_{k_x +i } } ^{j-i} $, and hence
$\Psi \left( \Phi ^{k_x +i } (l)\right)
=
\Phi |_{B_{k_x +i } } ^{j-i}
\left( \Phi ^{k_x +i } (l)\right)
=
\Phi ^{k_x +j } (l)$.
Thus
$\Psi $ maps
$\left\{ \Phi ^{k_x} (l) , \ldots , \Phi ^{{k_{x+1}}-1} (l)\right\} $
to itself.

Because
$\Psi $ maps
$P\setminus
\left\{ \Phi ^{v+k_x } (l), \ldots , \Phi ^{n-1} (l)\right\} $
to itself,
and because
$\Psi $ maps
$\left\{ \Phi ^{k_x} (l) , \ldots , \Phi ^{{k_{x+1}}-1} (l)\right\} $
to itself,
$\Psi $ must map the elements of
$\left\{ \Phi ^{v+k_x } (l), \ldots , \Phi ^{n-1} (l)\right\}
\setminus \{ a\} $
that are above elements of
$\left\{ \Phi ^{k_x} (l) , \ldots , \Phi ^{{k_{x+1}}-1} (l)\right\} $
to the elements of
$\left\{ \Phi ^{v+k_x } (l), \ldots , \Phi ^{n-1} (l)\right\}
\setminus \{ b\} $
that are above elements of
$\left\{ \Phi ^{k_x} (l) , \ldots , \Phi ^{{k_{x+1}}-1} (l)\right\} $.
By Lemma \ref{abovel}, no element of
$\left\{ \Phi ^{v+k_{x+1} } (l), \ldots , \Phi ^{n-1} (l)\right\} $
is above an element of
$\left\{ \Phi ^{k_x} (l) , \ldots , \Phi ^{{k_{x+1}}-1} (l)\right\} $.
Because
$k_{x+1} \leq \min \{ i_a , i_b \} $,
neither of $a$ or $b$ is in
$\left\{ \Phi ^{v+k_x} (l) , \ldots , \Phi ^{v+{k_{x+1}}-1} (l)\right\} $,
and hence $\Psi $ maps
$\left\{ \Phi ^{v+k_x} (l) , \ldots , \Phi ^{v+{k_{x+1}}-1} (l)\right\} $
to itself.

Now
let
$k_{x} \leq \min \{ i_a , i_b \}
<k_{x+1}
$.
Because
the {\em Claim}
holds for all
$z\in \{ -1, \ldots , x-1\} $, the same proof as
above
shows that
$\Psi $ maps
$\left\{ \Phi ^{k_x} (l) , \ldots , \Phi ^{{k_{x+1}}-1} (l)\right\} $
to itself.
Because $d_a $ is the only lower bound of
$a$ with $|\uparrow l|$ upper bounds in $P$,
$d_a =\Phi ^{i_a } (l)$
is the unique element of
$\left\{ \Phi ^{k_x} (l) , \ldots , \Phi ^{{k_{x+1}}-1} (l)\right\} $
with the fewest upper bounds in
$P\setminus \{ a\} $.
Because $\Psi $ is an isomorphism,
$\Psi $ must map
$d_a $ to the unique element of
$\left\{ \Phi ^{k_x} (l) , \ldots , \Phi ^{{k_{x+1}}-1} (l)\right\} $
with the fewest upper bounds in
$P\setminus \{ b\} $.
Because
all these elements have
$|\uparrow l|$ upper bounds in $P$,
and because
$d_b $ is the only lower bound of
$b$ with $|\uparrow l|$ upper bounds in $P$,
$d_b $ must be in
$\left\{ \Phi ^{k_x} (l) , \ldots , \Phi ^{{k_{x+1}}-1} (l)\right\} $
and
the $\Psi $-image of $d_a $ must be $d_b $.
\qed

\vspace{.1in}

{\bf Proof of Theorem \ref{minmaxpsrecon}.}
By
Lemma \ref{partitionPsi},
an ordered set $P$ as in
Lemma \ref{dnewisolem} cannot exist.
\qed

\section{Conclusion}
\label{conclusionsec}

Theorem \ref{minmaxpsrecon} is another entry in the
ever-growing list of reconstructible ordered sets.
For the task of reconstructing the minimal and maximal decks,
Theorem \ref{minmaxpsrecon} eliminates any and all concerns regarding
isomorphisms between minimal and maximal cards.
Theorem \ref{getnonexcndeck},
Corollary \ref{getsomemincards},
and Lemma \ref{identnonbotofSUMRC}
are further steps for this task:
Theorem \ref{getnonexcndeck} identifies
the deck that consists of the maximal and the minimal cards
together and
Corollary \ref{getsomemincards} and
Lemma \ref{identnonbotofSUMRC} identify
some maximal
cards.
We obtain the following
consequence, which, in light of
Theorem 6.6
of \cite{Schmarked}
(ordered sets of width $3$ are reconstructible if their
maximal cards on which $P\setminus \max (P)$ is ``marked"
are reconstructible),
is a significant
step towards finishing the
reconstruction of ordered sets of width $3$,

\begin{cor}
\label{width3mindeck}

Let $P$ be an ordered set of width $3$.
Then the (unmarked) minimal and maximal decks of $P$
are reconstructible.

\end{cor}

{\bf Proof.}
By Theorem 7.4 in \cite{Schmarked}, ordered sets
of width $3$ with $2$ maximal elements are
reconstructible.
Hence, by duality,
the only case we need to consider
is
that
$P$ has $3$ minimal elements and $3$ maximal elements.
Moreover, because by \cite{JHthes}
ordered sets with $\leq 11$
elements are reconstructible,
we can assume that $|P|\geq 12$.
Finally, by Corollary \ref{getsomemincards}, we only need to
consider the case that
$P\setminus \min (P) $ has a component
that has a
minmax pair $(l,h)$ of pseudo-similar points
or that is a singleton $\{ l\} $, and
$l$ has a unique lower cover.
We shall freely use duality.
Note that
$P\setminus \min (P)$
is reconstructible by Lemma \ref{smallcardscor}.

{\em Case 1.
$P\setminus \min (P)$
is connected and has a minmax pair of pseudo-similar
points.}
In this case, $P$ is reconstructible by Lemma \ref{getP-lLemma}.

{\em Case 2.
$P\setminus \min (P)$
is disconnected and
has a component $Q$ that is not a chain
and
that has a minmax pair $(l,h)$
of pseudo-similar
points.}
In this case,
$P\setminus \min (P)$
consists of $Q$ and a chain $C$ and
both sets (clearly) are identifiable on
$P\setminus \min (P)$.
By Theorem \ref{getnonexcndeck}, we can identify the
``rank $1$ deck" ${\cal N} _P ^1 $
of nonmaximal cards whose removed point has rank $1$.

{\em Case 2.1:
$|C|\geq 2$.}
Let $y$ be the smallest element of $C$.
In ${\cal N} _P ^1 $, we can
identify the card
$P\setminus \{ y\} $ as the card $Y$ in ${\cal N} _P ^1 $
such that
$Y\setminus \min (Y)$ consists of
a component isomorphic to $Q$ and a chain with
one less element than $C$.
Because $Q$ is rigid,
$(P\setminus \min (P))\setminus \{ y\} $ is rigid,
Because $y$ is identifiable on
any card
$Z:=P\setminus \{ z\} $ as in Lemma \ref{smallcardscor},
$P$ is then
reconstructible via Lemma \ref{identcardinupperlevels}.

{\em Case 2.2: $Q\setminus \{ l\} $ is connected.}
In ${\cal N} _P ^1 $, we can
identify the card
$P\setminus \{ l\} $ as the card $Y$ in ${\cal N} _P ^1 $
such that
$Y\setminus \min (Y)$ consists of
a rigid set isomorphic to $Q\setminus \{ l\} $
(which is not a chain)
and a chain with
as many elements as $C$.
Because $l$ is identifiable on
any card
$Z:=P\setminus \{ z\} $ as in Lemma \ref{smallcardscor}
and
$(P\setminus \min (P))\setminus \{ l\} $ is rigid,
$P$ is then
reconstructible via Lemma \ref{identcardinupperlevels}.

{\em Case 2.3: $|C|=1$ and $Q\setminus \{ l\} $ is disconnected.}
When
$Q\setminus \{ l\} $ is disconnected,
because $Q$ has width $2$, $Q$ is the
disjoint union of two chains, which are $K_h $ and $L$.
In this case,
$l<L$, $L\cup \{ l\} =K_l $
and
$l$ is a lower cover of $h$.
Let $s$ be the smallest element of $K_h $.
Because $|P|\geq 12$, we have that $|Q|\geq 8$ and
hence $|K_l |=|K_h |\geq 4$,
which means that $Q\setminus \{ s\} $ is rigid.
In ${\cal N} _P ^1 $, we can
identify the card
$P\setminus \{ s\} $ as the card $Y$ in ${\cal N} _P ^1 $
such that
$Y\setminus \min (Y)$ consists of
a set isomorphic to $Q\setminus \{ s\} $
and a singleton component.
Because $s$ is identifiable on
any card
$Z:=P\setminus \{ z\} $ as in Lemma \ref{smallcardscor}
and
$(P\setminus \min (P))\setminus \{ s\} $ is rigid,
$P$ is then
reconstructible via Lemma \ref{identcardinupperlevels}.
The completes {\em Case 2}.

{\em Case 3.
$P$ has a ranging chain
$C_b $, call its only minimal lower bound $b$,
and a dually ranging chain
$C_t $, call its only maximal upper bound $t$.}
Note that,
by {\em Cases 1 and 2}, by
Corollary \ref{getsomemincards},
and by their respective duals, this is
the only remaining case in which the minimal and maximal decks
cannot be identified.

Let $u$ be the top element of $C_b $
and let $d$ be the bottom element of $C_t $.
Let $M$ be the set of all elements that are not comparable to
any element in $C_b \cup C_t $.
Because
$P$ is connected,
we claim that $u$ and $d$ are not comparable to
each other:
Indeed $d<u$ would imply that $d$ is the only lower bound of
$C_b $, that $u$ is the only upper bound of $C_t $,
and then the chain $C_b \cup C_t $ would be a component of $P$.

Because $P$ is connected of width $3$
and because $u$ and $d$ are not comparable to
each other or to any element of $M$,
we conclude that $M$ is a chain.
Let $a$ be the third minimal element of
$P$ and note that, because $C_b $ is ranging with
unique strict lower bound $b$,
we have $a\in M$.
Dually, let
$v$ be the
third maximal element of
$P$ and note that $v\in M$.
Because $P=C_b \cup M\cup C_t \cup \{ b,t\} $
and $|P|\geq 12$,
we obtain that
$|C_b |+|M|+|C_t|\geq 10$

From
$|C_b |+|M|+|C_t|\geq 10$, we obtain
$|M|\geq 6$ or $|C_b |\geq 2$ or $|C_t |\geq 2$.
Clearly, if $|C_b |\geq 2$ or $|C_t |\geq 2$,
then $P$ is decomposable.
Note that the only elements of $P$ that
could be
comparable to
elements of $M$ are $b$ and $t$.
Consequently,
in case $|M|\geq 6$,
$M$ contains a nontrivial order-autonomous chain
$K$ that is a ranging chain with $\geq 5$ elements, a
dually ranging chain with $\geq 5$ elements,
or $K$ is neither ranging nor dually ranging with $\geq 2$ elements.
In particular,
in case $|M|\geq 6$,
$P$ is decomposable, too.

Because the number of (dually) ranging chains
does not change when an element of a nontrivial order-autonomous set is
removed,
the number of (dually) ranging chains
can be obtained from an NTMA-card.
All sizes of order-autonomous
chains can be identified by
part \ref{decomplem5} of Lemma \ref{decomplem}
(and in case 2 elements are in nontrivial order-autonomous subsets, the size is $2$).
In particular, independent of whether we have 2 or 3 ranging or
dually ranging chains, we can
identify the sizes of the ranging and the dually ranging chains,
which then means that all cases below are identifiable.

{\em Case 3.1: $P$ has exactly one ranging
chain $C_b $ and $|C_b |\geq 2$ or $P$ has exactly one
dually ranging
chain $C_t $ and $|C_t |\geq 2$.}
By duality, we only need to consider the case that
$C_b $ is the unique ranging
chain in $P$ and $|C_b |\geq 2$.
Let $P\setminus \{ x\} $ be an NTMA-card with
as few elements in its unique ranging chain
as possible.
Then
$x\in C_b $ and
we can reconstruct $P$
by, on $P\setminus \{ x\} $, adding a new top element to
the identifiable order-autonomous set
$C_b \setminus \{ x\} $.

{\em Case 3.2.
$P$ has two ranging chains or two
dually ranging chains.}
By duality, we only need to consider the case that
$P$ has two ranging chains.
In this case,
there is a unique dually ranging chain $C_t $
and, by the above,
we can assume $|C_t |=1$.
Consequently,
one of the
ranging chains
has at least $4$ elements.
Let the ranging chains be $K_1 $ and $K_2 $ with
$|K_1 |\geq |K_2 |$.
In the identifiable case $|K_2 |\geq 2$,
we reconstruct $P$ from an identifiable
NTMA-card $P\setminus \{ x\} $ with a
ranging chain with $|K_1 |$ elements
and another
ranging chain with $|K_2 |-1$ elements
by realizing that the latter chain
must be $K_2 \setminus \{ x\} $ and by adding a new top element to it.
In the identifiable case $|K_2 |=1$, we have $|K_1 |\geq 7$
and
we reconstruct $P$ from an identifiable
NTMA-card $P\setminus \{ x\} $ with a
ranging chain with $|K_1 |-1$ elements
and another
ranging chain with $1$ element
by realizing that the former chain
must be $K_1 \setminus \{ x\} $ and by adding a new top element to it.

{\em Case 3.3.
$P$ has exactly one ranging chain $C_b $,
exactly one dually ranging chain $C_t $,
and $|C_b |=|C_t |=1$.}
In this final case, $|M|\geq 8$, and
$M$ contains at least one order-autonomous chain
with at least $3$ elements.
Moreover, because any automorphism must map
(dually) ranging chains to
(dually) ranging chains,
every automorphism of the index set $T$ of $P$
must map $I[M]$ to itself.
Consequently,
the index set of the canonical decomposition
of $P$ is rigid, and
$P$ is reconstructible
via
part \ref{decomplem6} of Lemma \ref{decomplem}.
\qed

\vspace{.1in}

We also have a new recognition result as follows.
Recall that $x\in P$ is called {\bf irreducible} iff
$x$ has a unique upper cover or a unique lower cover.
$P$ is called {\bf dismantlable} iff there is an
enumeration
$P=\{ x_1 , \ldots , x_n \} $ such that, for all
$j\in \{ 1, \ldots , n-1\} $, the element
$x_j $ is irreducible in $P\setminus \{ x_1 , \ldots , x_{j-1} \} $.

\begin{cor}

Dismantlable ordered sets are recognizable.

\end{cor}

{\bf Proof.}
By
part \ref{nhoodfacts3} of Lemma \ref{nhoodfacts},
the deck of maximal cards $P\setminus \{ x\} $
such that $x$ has a unique
lower
cover is
reconstructible.
If any such card is dismantlable, then $P$ is dismantlable.
Dually, the deck of
minimal cards $P\setminus \{ x\} $
such that $x$ has a unique
upper
cover is
reconstructible, and,
if any such card is dismantlable, then $P$ is dismantlable.

By
part \ref{nhoodfacts1} of Lemma \ref{nhoodfacts},
we can reconstruct the isomorphism types of the
neighborhoods of all points of rank $k$ and their multiplicities.
By
Theorem \ref{getnonexcndeck}, we can determine the
decks of nonextremal cards obtained by
removal of a nonextremal element of rank $k$.
For each card
$P\setminus \{ x\} $ in this ``nonextremal rank $k$ deck,"
the isomorphism type of the
neighborhood
of $x$
is the unique isomorphism type
$N$ of a neighborhood
for a point of rank $k$ in $P$ such that
$P\setminus \{ x\} $
has fewer points of rank $k$
whose neighborhoods
are isomorphic to $N$ than $P$ does.
In particular,
we can determine
if $x$ had a unique upper cover or a
unique lower cover.
If any nonextremal card
such that $x$ had a unique upper or lower cover
is dismantlable, then $P$ is dismantlable.
If none of the above holds, then $P$ is not dismantlable.
\qed

\vspace{.1in}

If we could reconstruct the maximal and minimal decks, then,
similar to the preceding proof,
we could reconstruct the isomorphism type of the neighborhood
of the removed point for every card.
In particular, we would then know
the number of upper and lower bounds for the removed point.
This would then make the reconstruction conjecture for ordered sets a
special case of the
reconstruction conjecture in Section 3 of
\cite{Ramadacard}, which
poses that all digraphs should be reconstructible from
their decks provided that, for each card, the
in-degree and the out-degree of the removed point is known, too.
Indeed, for order reconstruction, availability of the
isomorphism type of the neighborhood would
provide even more information.

\section{Declarations}

\begin{itemize}
\item
{\bf Ethical Approval.}
Not applicable, because no human or animal studies were involved.

\item
{\bf Competing Interests.}
The author has no relevant financial or non-financial
interests to disclose.

\item
{\bf Authors' Contributions.}
This is a single author paper with all
work done by the author.
Contributions by non-authors are acknowledged above.

\item
{\bf Funding.}
No funds, grants, or other support was received.

\item
{\bf Availability of Data and Materials.}
Data sharing is not applicable, as no datasets
were generated or analysed herein.

\end{itemize}

{\bf Acknowledgement.}
I thank an anonymous referee for not accepting an
earlier version
of the proof of
Theorem \ref{getnonexcndeck}
that was almost twice as long and looked like
``spaghetti code."
The paper is better for it.


\begin{thebibliography}{99}













\bibdiss{JHthes}
{J. Hughes}{2004}{The Computation and Comparison of Decks of Small Ordered Sets}
{MS thesis}{Louisiana Tech University}



\bibcoll{Illedec}
{P. Ille}
{1993}
{Recognition problem in reconstruction for decomposable
relations}{N. W. Sauer et. al. (eds.)}{Finite and
Infinite Combinatorics in Sets and Logic}{Kluwer Academic
Publishers}{189--198}





\bibart{KRtow}
{D. Kratsch and J.-X. Rampon}{1994}
{Towards the reconstruction of
posets}{Order}{11}{317--341}



\bibart{Ramadacard}
{S.Ramachandran}{1981}{On a new digraph reconstruction conjecture}{J. Combin. Theory Ser. B}
{31}{143--149}





\bibart{JXsurvey}
{J.-X. Rampon}{2005}{What is reconstruction for ordered sets?}
{Discrete Math.}{291}{191--233}




\bibart{Schneighdec}
{B. Schr\"oder}{2000}
{Reconstruction of the neighborhood deck of ordered sets}
{Order}{17}{255--269}










\bibart{Schmarked}
{B. Schr\"oder}{2003}
{On Ordered Sets with Isomorphic Marked Maximal Cards}
{Order}
{20}{299--327}





\bibart{bman1}
{B. Schr\"oder}{2010}
{Pseudo-Similar Points in Ordered Sets}
{Discrete Mathematics}{310}{2815--2823}








\bibart{SchrSetRec}
{B. Schr\"oder}{2022}
{Set recognition of decomposable graphs and steps
towards their reconstruction}
{Abh. Math. Semin. Univ. Hambg.}{92}{1--25}

\url{https://doi.org/10.1007/s12188-021-00252-0}



\end{thebibliography}
\end{document}